%% file: simplex3d.tex
\documentclass[12pt]{article}
\date{}

\usepackage{anysize,amsthm,amsmath,amsfonts,amssymb,graphicx,calc,url,overpic,titling,paralist,color,epsfig}
\graphicspath{{EPS/}{EPS-A/}}

\theoremstyle{plain}
\newtheorem{lemma}{Lemma}[section]
\newtheorem{theorem}[lemma]{Theorem}

\newtheorem{proposition}[lemma]{Proposition}
\theoremstyle{definition}

\theoremstyle{remark}
 
\newcommand{\R}{\mathbb{R}}
\newcommand {\st} {\text{start}}
\newcommand {\GD} {\text{GD}}
\newcommand {\RE} {\text{RE}}

\newcommand {\incut} {\operatorname{\delta}^{\text{in}}\!}
\newcommand {\outcut} {\operatorname{\delta}^{\text{out}}\!}
\newcommand\EXP{\mathbf{E}}

\thanksmarkseries{arabic}
\title{The Simplex Algorithm in Dimension Three%
  \thanks{%
    Work on this paper by Micha Sharir was supported by NSF Grants
    CCR-97-32101 and CCR-00-98246, by a grant from the U.S.-Israeli
    Binational Science Foundation, by a grant from the Israel Science
    Fund (for a Center of Excellence in Geometric Computing), and by
    the Hermann Minkowski--MINERVA Center for Geometry at Tel Aviv
    University.  
    Volker Kaibel and Rafael Mechtel were supported by the 
    DFG-Forschergruppe \emph{Algorithmen, Struktur, Zufall} 
    (FOR 413/1-1, Zi 475/3-1).  
    G\"unter M. Ziegler acknowledges partial support by the
    Deutsche Forschungs-Gemeinschaft (DFG), FZT86, ZI 475/3 and ZI 475/4
    and by the GIF project 
    \emph{Combinatorics of Polytopes in Euclidean Spaces} (I-624-35.6/1999)
    Part of the work was done during the workshop 
   ``Towards the Peak'' at La Claustra, Switzerland, August 2001.
  }%
}
\author{%
  Volker Kaibel%
  \thanks{%
    DFG Research Center ``Mathematics for Key Technologies'', 
    MA 6--2, TU Berlin, 10623~Berlin, Germany;
    \texttt{kaibel@math.tu-berlin.de}
  }
\and 
  Rafael Mechtel%
  \thanks{%
    MA 6--2, TU Berlin, 10623~Berlin, Germany; 
    \texttt{\{mechtel,ziegler\}@math.tu-berlin.de}
  }
\and
  Micha Sharir\thanks{%
    School of Computer Science, Tel Aviv University, Tel-Aviv
    69978, Israel and Courant Institute of Mathematical Sciences, 
    New York University, New York, NY 10012, USA; 
    \texttt{michas@post.tau.ac.il} 
  }
\and 
  G\"unter M. Ziegler\thanksmark{3}%
}

\begin{document}
\maketitle

\input intro

\input diam

\input rand3d

\input other


\paragraph{Acknowledgements.}
We are grateful to Emo Welzl and G\"unter Rote for inspiring
discussions and helpful comments.


\bibliographystyle{siam}


\end{document}

%% file: intro.tex
\begin{abstract}
  We investigate the worst-case behavior of the simplex algorithm 
  on linear programs with three variables, that is,
  on $3$-dimensional simple polytopes.
  Among the pivot rules that we consider, the ``random edge''
  rule yields the best asymptotic behavior as well as the most
  complicated analysis.
  All other rules turn out to be much easier to study, but
  also produce worse results: Most of them show essentially
  worst-possible behavior; this includes both Kalai's
  ``random-facet'' rule, which 
  without dimension restriction is known to be subexponential, as well as 
  Zadeh's deterministic history-dependent rule,
  for which no non-polynomial instances
  in  general dimensions have been found so far. 
\end{abstract}


\section{Introduction}
\label{sec:intro}

The simplex algorithm is a fascinating method for at least three
reasons: For computational purposes it is 
still the most efficient general
tool for solving linear programs, from a complexity point of view it
is the most promising candidate for a strongly polynomial time linear
programming algorithm, and last but not least, geometers are pleased
by its inherent use of the structure of convex polytopes.

The essence of the method can be described geometrically: Given a
convex polytope~$P$ by means of inequalities, a linear
functional~$\varphi$ ``in general position,'' and some
vertex~$v_{\text{start}}$, the simplex algorithm
chooses an edge to a neighboring vertex
along which $\varphi$ decreases strictly. Iterating this yields a
$\varphi$-monotone edge-path.  Such a path can never get stuck, and
will end at the unique $\varphi$-minimal (``optimal'') vertex of~$P$.

Besides implementational challenges, a
crucial question with respect to efficiency asks for a suitable
\emph{pivot rule} that prescribes
how to proceed with the monotone path at any vertex. 
Since Dantzig invented the simplex algorithm in the late
1940's \cite{Dan63}, a great variety of pivot rules have
been proposed. Most
of them (including Dantzig's original ``largest coefficient rule'')
have subsequently been shown to lead to exponentially long paths in
the worst case. (See~\cite{AZ96} for a survey.)  Prominent
exceptions are Zadeh's history-dependent ``least entered'' rule, and 
several randomized pivot rules. Particularly remarkable is the
``random facet'' rule proposed by Kalai~\cite{Kal92}; its 
expected path length for all instances
is bounded subexponentially in the number of facets.
See also Matou\v{s}ek et al.~\cite{MSW96}.

In this paper, we analyze the worst-case behavior of
the simplex method on $3$-\allowbreak dimen\-sional 
simple polytopes for some
well-known pivot rules. 
At first glance, the $3$-dimensional case may seem trivial,
since by Euler's formula
a $3$-polytope with $n$ facets has at most $2n-4$ vertices
(with equality if and only if the polytope is simple),
and there are examples where $n-3$ steps are needed
for any monotone path to the optimum (see, e.g., Figure~\ref{fig:li}).
Therefore, for \emph{any} pivot rule the simplex algorithm
is linear, with at least $n-3$ and at most $2n-5$ steps in the worst case.
However, no pivot rule is known that would work with
at most $n-3$ steps.


In order to summarize our results, we define the following measure of
quality. Fix a pivot rule~$\mathcal{R}$. For every
$3$-dimensional polytope $P\subset\R^3$ and for every linear
functional $\varphi:\R^3\longrightarrow\R$ \emph{in general position}
with respect to $P$
(i.\,e., no two vertices of~$P$ have the same $\varphi$-value), denote
by $\lambda_{\mathcal{R}}(P,v_{\text{start}})$ the path length
(expected path length, if~$\mathcal{R}$ is randomized) produced by the
simplex algorithm with the pivot rule~$\mathcal{R}$, when started at
vertex~$v_{\text{start}}$. The \emph{linearity coefficient}
of~$\mathcal{R}$ is
\[
\Lambda(\mathcal{R})\ :=\ \limsup_{n(P)\rightarrow\infty}
  \Big\{
    \frac{\lambda_{\mathcal{R}}(P,v_{\text{start}})}{n(P)}\ :\
    P,\varphi,v_{\text{start}}\text{ as above}
  \Big\},
\]
where $n(P)$ is the number of facets of~$P$. With the usual
simplifications
for a geometric analysis (cf.  \cite{KlKl1}, \cite[Lect.~3]{Zie95},
\cite{AZ96}),
we may restrict our attention to \emph{simple} $3$-dimensional
polytopes~$P$ (where each vertex is contained in precisely $3$
facets). So we only consider $3$-dimensional polytopes $P$,
with~$n=n(P)$ facets, $3n-6$ edges, and $2n-4$ vertices.
By the discussion above, the linearity coefficient satisfies
$1\le\Lambda(\mathcal{R})\le 2$ for every pivot rule~$\mathcal{R}$.

The most remarkable aspect of the picture that we obtain, in
Section~\ref{sec:rand3d}, is that the ``random edge'' rule (``RE'' for
short) performs quite well (as it is conjectured for general dimensions),
but it is quite tedious to analyze (as it has already been observed
for general dimensions).  The following bounds
for the random edge rule
\[
1.3473 \ \ \le \ \ \Lambda(\RE)\ \ \le\ \ 1.4943
\]
are our main results.
Thus we manage to separate $\Lambda(\RE)$ from the rather easily
achieved lower bound of $\frac43$, as well as from the already
non-trivial upper bound of $\frac32$.

On the other hand, in Section~\ref{sec:other}
we prove that the linearity coefficient
for the ``greatest decrease'' pivot rule is
$\Lambda(\GD)=\frac32$, while many other
well-known rules have linearity coefficient
$\Lambda=2$, including  the
    largest coefficient, 
    least index, 
    steepest decrease, and the 
    shadow vertex rules, as well as Zadeh's history-dependent
    least entered rule (not known to be super-polynomial
    in general),
    and Kalai's random facet rule (known to be sub-exponential in general).


%% file: diam.tex

\section{Basics}
\label{sec:monoHirsch3d}

Klee~\cite{Kle65} proved in 1965 that the ``monotone
Hirsch conjecture'' is true for $3$-dimensional polytopes, that is,
whenever the graph of a $3$-dimensional polytope~$P$ with~$n$ facets is
oriented by means of a linear functional in general position
there is 
a monotone path of length at most $n-3$ from any vertex to the sink~$v_{\min}$.
(See Klee \& Kleinschmidt \cite{KlKl1} for a survey of the 
Hirsch conjecture and its ramifications.)
Unfortunately, Klee's proof is not based on a pivot rule. 


\begin{theorem}[Klee~\cite{Kle65}]
  \label{thm:nonRev3d}
  For any simple $3$-polytope~$P\subset\R^3$, a linear
  functional $\varphi:\R^3\longrightarrow\R$ in general position
  for~$P$, and any vertex~$v_\st$ of~$P$, there is a
  $\varphi$-monotone path from~$v_\st$ to the $\varphi$-minimal vertex
  $v_{\min}$ of~$P$ that does not revisit any facet.
  
  In particular, there is a $\varphi$-monotone path from~$v_\st$ to
  $v_{\min}$ of length at most $n-3$.
\end{theorem}

It is not too hard to come up with examples showing that the bound
provided by Theorem~\ref{thm:nonRev3d} is best possible. One of
the constructions will be important for our treatment later on, so we
describe it below in Figure~\ref{fig:li}. 

A particularly useful tool for constructing
LP-oriented $3$-polytopes is the following result due to Mihalisin and
Klee.  It is stated in a slightly weaker version in their
paper, but their proof actually shows the following.

\begin{theorem}[Mihalisin \&\ Klee~\cite{MK00}]
  \label{thm:MK}
  Let $G=(V,E)$ be a planar $3$-connected graph,
  $f:V\longrightarrow\R$ any injective function, and denote by
  $\vec{G}$ the acyclic oriented graph obtained from~$G$ by directing
  each edge to its endnode with the smaller $f$-value.
  Then the following are equivalent:
\begin{compactenum}[\rm1.]
\item
  There exist a polytope $P\subset\R^3$ and a linear functional
  $\varphi:\R^3\longrightarrow\R$ in general position for~$P$, such
  that~$G$ is isomorphic to the graph of~$P$ and, for every $v\in V$,
  $f(v)$ agrees with the $\varphi$-value of the vertex of~$P$
  corresponding to~$v$.
\item Both {\rm(a)} and\/ {\rm(b)} hold:  
  \begin{compactenum}[\rm(a)]
  \item $\vec{G}$ has a unique sink in every facet (induced
    non-separating cycle) of~$G$, and
  \item there are three node-disjoint monotone paths joining
    the (unique) source to the (unique) sink of~$\vec{G}$.
  \end{compactenum}
\end{compactenum}
\end{theorem}

\noindent
Here the fact that the source and the sink of~$\vec{G}$ are unique (referred
to in condition~(b)) follows from~(a);
cf. Joswig et al.~\cite{JKK02}.
Equipped with Theorem~\ref{thm:MK}, one readily verifies that the
family of directed graphs indicated in Figure~\ref{fig:li} $(n\ge4)$
can be realized as convex $3$-polytopes, with associated linear functionals, demonstrating  that Klee's bound of $n-3$ on the length of a shortest
monotone path cannot be improved.

  \begin{figure}[ht]
    \begin{center}
      \input{EPS/li.pstex_t}
      \caption{ \label{fig:li}
        A worst case example for Klee's theorem, starting at $v_{n-3}$
        (and for Bland's rule, starting at $v_{2n-6}$; see
        Section~\ref{subsec:bland}).
        All edges are oriented from left to right.}
    \end{center}
  \end{figure}
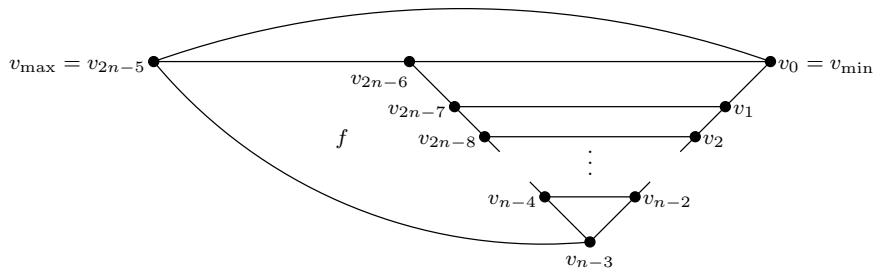

%% file: EPS/li.pstex_t
\begin{picture}(0,0)%
\includegraphics{li.pstex}%
\end{picture}%
\setlength{\unitlength}{4144sp}%
\begingroup\makeatletter\ifx\SetFigFont\undefined
\def\x#1#2#3#4#5#6#7\relax{\def\x{#1#2#3#4#5#6}}%
\expandafter\x\fmtname xxxxxx\relax \def\y{splain}%
\ifx\x\y   
\gdef\SetFigFont#1#2#3{%
  \ifnum #1<17\tiny\else \ifnum #1<20\small\else
  \ifnum #1<24\normalsize\else \ifnum #1<29\large\else
  \ifnum #1<34\Large\else \ifnum #1<41\LARGE\else
     \huge\fi\fi\fi\fi\fi\fi
  \csname #3\endcsname}%
\else
\gdef\SetFigFont#1#2#3{\begingroup
  \count@#1\relax \ifnum 25<\count@\count@25\fi
  \def\x{\endgroup\@setsize\SetFigFont{#2pt}}%
  \expandafter\x
    \csname \romannumeral\the\count@ pt\expandafter\endcsname
    \csname @\romannumeral\the\count@ pt\endcsname
  \csname #3\endcsname}%
\fi
\fi\endgroup
\begin{picture}(3780,1573)(721,-947)
\put(4501,254){\makebox(0,0)[lb]{\smash{\SetFigFont{8}{9.6}{rm}{\color[rgb]{0,0,0}$v_0 = v_{\text{min}}$}%
}}}
\put(3376,-376){\makebox(0,0)[b]{\smash{\SetFigFont{8}{9.6}{rm}{\color[rgb]{0,0,0}$\vdots$}%
}}}
\put(2521,-16){\makebox(0,0)[rb]{\smash{\SetFigFont{8}{9.6}{rm}{\color[rgb]{0,0,0}$v_{2n-7}$}%
}}}
\put(2701,-196){\makebox(0,0)[rb]{\smash{\SetFigFont{8}{9.6}{rm}{\color[rgb]{0,0,0}$v_{2n-8}$}%
}}}
\put(3061,-556){\makebox(0,0)[rb]{\smash{\SetFigFont{8}{9.6}{rm}{\color[rgb]{0,0,0}$v_{n-4}$}%
}}}
\put(3691,-556){\makebox(0,0)[lb]{\smash{\SetFigFont{8}{9.6}{rm}{\color[rgb]{0,0,0}$v_{n-2}$}%
}}}
\put(4051,-196){\makebox(0,0)[lb]{\smash{\SetFigFont{8}{9.6}{rm}{\color[rgb]{0,0,0}$v_{2}$}%
}}}
\put(4231,-16){\makebox(0,0)[lb]{\smash{\SetFigFont{8}{9.6}{rm}{\color[rgb]{0,0,0}$v_{1}$}%
}}}
\put(1891,-196){\makebox(0,0)[b]{\smash{\SetFigFont{8}{9.6}{rm}{\color[rgb]{0,0,0}$f$}%
}}}
\put(721,254){\makebox(0,0)[rb]{\smash{\SetFigFont{8}{9.6}{rm}{\color[rgb]{0,0,0}$v_{\text{max}} = v_{2n-5}$}%
}}}
\put(3376,-916){\makebox(0,0)[b]{\smash{\SetFigFont{8}{9.6}{rm}{\color[rgb]{0,0,0}$v_{n-3}$}%
}}}
\put(2296,164){\makebox(0,0)[rb]{\smash{\SetFigFont{8}{9.6}{rm}{\color[rgb]{0,0,0}$v_{2n-6}$}%
}}}
\end{picture}

%% file: rand3d.tex
\section{The Random Edge Rule} 
\label{sec:rand3d}

At any non-optimal vertex, the \emph{random edge} pivot rule 
takes a step to
one of its improving neighbors, chosen uniformly at random.
Thus the expected number $\EXP(v)$ of steps that the random edge rule would
take from a given vertex $v$ to the optimal one $v_{\min}$ may
be computed recursively as
\begin{equation} \label{eq:re-rec}
\EXP(v)\ =\ 1+\frac1{|\outcut(v)|}\sum_{u:(v,u)\in\outcut(v)} \EXP(u),
\end{equation} 
where $\outcut(v)$ denotes the set of edges that leave $v$
(that is, lead to \emph{better} vertices), so that
$|\outcut(v)|$ is the number of neighbors of $v$ whose $\varphi$-value
is smaller than that of~$v$.

Despite its simplicity and its (deceptively) simple recursion, 
this rule has by now resisted several attempts to
analyze its worst-case behavior, with a few exceptions for special
cases, namely linear assignment problems (Tovey~\cite{Tov86}), the
Klee-Minty cubes (Kelly~\cite{Kel81}, 
G\"artner et al.~\cite{GHZ98}), and $d$-dimensional linear programs with at
most $d+2$ inequalities (G\"artner et al.~\cite{GSTVW01}). All known
results leave open the possibility that the expected number of steps
taken by the random edge rule  on a $d$-dimensional linear program with $n$ inequalities  could be bounded by a polynomial,
perhaps even by $O(n^2)$ or~$O(dn)$, 
where $n$ is the number of facets. 

However, Matou\v{s}ek and Szabo~\cite{MS04} recently showed that the random edge rule does not have a polynomially bounded running time on the larger class of \emph{acyclic unique sink orientations} (\emph{AUSO}'s), i.e., acyclic orientations of the graph of a polytope that induce unique sinks in all non-empty faces (cf. condition 2(a) in Theorem~\ref{thm:MK}). They exhibited particular AUSO's on $d$-dimensional cubes for which random edge needs at least $\text{const}\cdot 2^{\text{const}\cdot d^{1/3}}$ steps.

\subsection{Lower Bounds}


The lower bound calculations appear to be much simpler if we do not
use the recursion given above, but instead use a ``flow model.'' For
this, fix a starting vertex $v_{\text{start}}$, and denote by
$p(v)$ the probability that the vertex $v$ will be visited by a random
edge path from $v_{\text{start}}$ to $v_{\min}$, and similarly by
$p(e)$ the probability that a directed edge $e$ will be traversed.
Then the probability that a vertex $v$ is visited is
the sum of the probabilities that the edges leading into $v$ are traversed,
\[
p(v)\ =\ \sum_{e\in\incut(v)} p(e)
\]
if $v$ is not the starting vertex. 
(Here $\incut(v)$ denotes the set of edges that enter~$v$.)
Furthermore, by definition
of the random edge rule we have 
\begin{equation} \label{eq:re-probs}
p(e)\ =\ \frac{1}{|\outcut(v)|} p(v)
\qquad\text{for all } e \in \outcut(v)
\end{equation}
at each non-optimal vertex.
The random edge rule thus induces a flow $\big(p(e)\big)_{e\in E}$ of
value~$1$ from $v_{\text{start}}$ to~$v_{\min}$.  The expected
path length $\EXP(v_{\text{start}})$ is then given by
\begin{equation} \label{eq:re-epl}
\EXP(v_{\text{start}})\ =\ \sum_{e \in E} p(e),
\end{equation}
and we refer to it as the \emph{cost} of the flow~$\big(p(e)\big)_{e\in E}$.

\begin{theorem}\label{thm:rand3dLowerBound}
The linearity coefficient of the random edge rule satisfies
\[
\Lambda({\rm RE})\ \ge\ \ \tfrac {1897}{1408}\ >\ 1.3473.
\]
\end{theorem}

\begin{proof}
We describe a family of LPs which show the above lower bound on the
linearity coefficient. We start with the graph of the dual-cyclic
polytope $C_3(k)^\Delta$ with the orientation depicted in
Figure \ref{fig:re-bb}, and refer to this as the \emph{backbone} of the
construction.

  \begin{figure}
    \begin{center}
      \input{EPS/backbone.pstex_t}
      \caption{ \label{fig:re-bb}
        Lower bound construction for the random edge rule: The
        backbone polytope.  All edges are oriented from left to right.
      }
    \end{center}
  \end{figure}
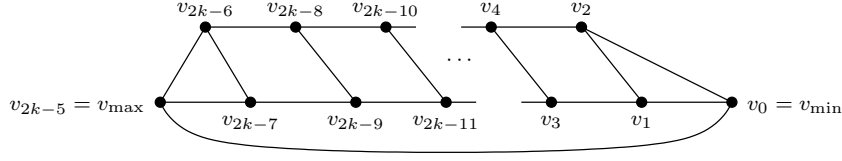

Starting at the vertex $v_{2k-7}$, the simplex algorithm 
will take the path along
the $k-2$ vertices $v_{2k-7},v_{2k-9},\dots,v_3,v_1,v_0$.
Replacing each vertex in the path
by a copy of the digraph depicted in Figure~\ref{fig:re-conf} 
--- called a \emph{configuration} in the following --- 
yields the desired LP. 
The corresponding feasible polytope can be constructed explicitly by applying 10 suitable successive
vertex cuts at each vertex~$v_i$ of the backbone. 
Alternatively, one can check that the orientations we get satisfy
the conditions of Theorem \ref{thm:MK}. 

  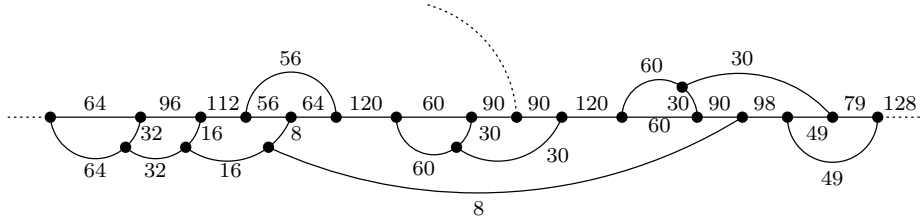
\begin{figure}[ht]
    \begin{center}
      \input{EPS/re.pstex_t}
      \caption{\label{fig:re-conf}%
        Lower bound construction for the random edge rule: The
        configuration.  All edges are oriented from left to right.
        The target of the rightmost edge (not shown) is the starting
        node of the next configuration. The middle dotted edge enters
        the configuration from the corresponding vertex of the top
        backbone row. The actual flow at each edge is $1/128$ times the number written next to the edge.}
    \end{center}
  \end{figure}

The maximal and minimal vertex of each configuration are visited with
probability~$1$. 
We send 128 units of flow (each of value $\frac1{128}$)
through each configuration 
according to~\eqref{eq:re-probs}; see Figure~\ref{fig:re-conf}. This
yields the flow-cost of 
$\tfrac {1897}{128}$ for each of the $k-2$ configurations. (The last
configuration produces flow-cost of $\tfrac {1897}{128} - 1$ only, as
it does not have a leaving edge.)

We take the maximal vertex of the configuration at $v_{2k-7}$
as the starting vertex~$v_\st$. Using equation~\eqref{eq:re-epl}
we obtain for the expected cost~$\EXP(v_\st)$:
\begin{gather*}
\EXP(v_\st) = (k-2)\tfrac{1897}{128} - 1.
\end{gather*}
With $n = k + 10(k-2)$ this yields
\begin{gather*}
\EXP(v_\st) \ =\ \tfrac{n-2}{11}\cdot\tfrac{1897}{128} - 1 
\ =\ \tfrac{1897}{1408} n - \tfrac{5202}{1408}, 
\end{gather*}
which proves the lower bound.
\end{proof}

The configuration depicted in Figure \ref{fig:re-conf} was found by
complete enumeration of the
acyclic orientations satisfying condition
(a) of Theorem~\ref{thm:MK} (AUSOs) 
on 3-polytopes with $n \le 12$ facets.
In particular, our proof of Theorem~\ref{thm:rand3dLowerBound}
includes a worst-possible example for $n=12$.
We refer to Mechtel \cite{Mechtel03} for more details of the search procedure,
as well as for a detailed analysis of properties of
worst-case examples for the random edge rule.

\subsection{Upper Bounds}\label{subsec:re-upper}

\begin{theorem} \label{thm:re-upper}
The linearity coefficient of the random edge rule satisfies
\[
\Lambda({\rm RE})\ \le\ \frac{130}{87} < 1.4943.
\]
\end{theorem}


\begin{proof}

Consider any linear program on a
simple $3$-polytope with $n$ facets, 
with a linear objective function $\varphi$ in general position.
We will refer to the $\varphi$-value of a vertex as its ``height.''
A \emph{$1$-vertex} will denote a vertex with exactly one 
neighbor that is lower with respect to~$\varphi$.
Similarly, a \emph{$2$-vertex} has exactly $2$ lower neighbors.
Consequently, from any $1$-vertex the random edge rule proceeds deterministically
 to the unique improving neighbor, 
and from any $2$-vertex  it proceeds to one of the two
improving neighbors, each with probability~$\frac12$.

Basic counting yields that our LP has exactly ($n-3$) $1$-vertices and
($n-3$) $2$-vertices in addition to the unique maximal vertex $v_{\max}$ 
and the unique minimal vertex $v_{\min}$, which have $3$ and $0$ lower 
neighbors, respectively.
For the following, we also assume that the vertices are
sorted and labelled $v_{2n-5},\dots,v_1,v_0$
in decreasing order of their objective function values,
with $v_{\max}=v_{2n-5}$ and $v_{\min}=v_0$.

For any vertex $v$, let $N_1(v)$ (resp., $N_2(v)$) denote the number
of $1$-vertices (resp., $2$-vertices) that are not higher than $v$
(including $v$ itself).  Put $N(v)=N_1(v)+N_2(v)$.  For all vertices
$v$ other than the maximal one this is the number of vertices lower
than $v$, that is, $N(v_i)=i$ for all $i\neq2n-5$.

We will establish the following generic inequality:
\begin{equation}
  \label{eq:generic}
  \EXP(v)\ \le\ \alpha N_1(v)+ \beta N(v)
\end{equation}
Here,  $\alpha$ and $\beta$ are constants whose values will be fixed
later. 

The proof of~\eqref{eq:generic} will proceed by induction on $N(v)$. The inductive step will be subdivided into 24
distinct cases.  Each case depends on a linear inequality on $\alpha$ and
$\beta$ that, when satisfied, justifies the induction step  in that case.  Since our case analysis
is complete, we have a proof of \eqref{eq:generic} for any pair
$(\alpha,\beta)$ that satisfies all the 24 inequalities.  

Because we always have
$N_1(v),N_2(v) \le n-3$, we obtain 
\[
\EXP(v)\ \le\ \alpha N_1(v) + \beta(N_1(v)+N_2(v)) =
(\alpha+\beta) N_1(v) + \beta N_2(v)\ \le\ (\alpha+2\beta)(n-3)
\]
for $v\neq v_{\max}$. The single vertex $v_{\max}$ is irrelevant for the asymptotic considerations. Thus we minimize $\alpha+2\beta$
subject to the linear constraints posed by the various cases; this
leads to an LP in two variables with 24 constraints, whose optimal
solution is $(\alpha,\beta)=(\frac{46}{87},\frac{42}{87})$, of
value~$\frac{130}{87}<1.4943$. This yields the upper bound on $\Lambda(\text{RE})$ stated in the theorem.

We will now prove~\eqref{eq:generic} by induction on $N(v)$. 
The
base case $N(v)=0$ is obvious, since $v$ is the optimum in this case,
and $\EXP(v)=0$.  Suppose now that~\eqref{eq:generic}  holds
for all vertices lower than some vertex $v$.

By an appropriate unwinding of the recursion (\ref{eq:re-rec}),
we express $\EXP(v)$ in terms of the expected cost $\EXP(w_i)$ of certain
vertices $w_i$ that are reachable from $v$ via a few downward edges. The
general form of such a recursive expression will be
\[
\EXP(v)\ =\ c + \sum_{i=1}^k \lambda_i \EXP({w_i}),
\]
where $\lambda_i > 0$ for $i=1,\dots,k$, and $\sum_i\lambda_i = 1$.

Since we assume by induction that 
$\EXP({w_i}) \le \alpha N_1(w_i) + \beta N(w_i)$, 
for each $i$, it suffices to show that
\[
\sum_{i=1}^k \alpha \lambda_i \left(N_1(v)-N_1(w_i)\right) +
\sum_{i=1}^k \beta \lambda_i \left(N(v)-N(w_i)\right)\ \ge\ c.
\]
Write
\[
\Delta_1(w_i) := N_1(v)-N_1(w_i),\qquad
\Delta(w_i) := N(v)-N(w_i),
\]
for $i=1,\ldots,k$. (These terms are defined with
respect to the vertex~$v$ that is currently considered.) Here
$\Delta(w_i)$ is the \emph{distance} between~$v$ and~$w_i$, that is,
one plus the number of vertices between $v$ and $w_i$
in the numbering of the vertices $(v_{2n-5},\dots,v_0)$ detailed above.
Clearly $\Delta(w_i)\ge\Delta_1(w_i)$.

We thus need to show that for each vertex $v$,
\begin{equation} \label{ind}
\alpha \sum_{i=1}^k \lambda_i \Delta_1(w_i) +
\beta \sum_{i=1}^k \lambda_i \Delta(w_i) \ge c.
\end{equation}

At this point we start our case analysis.

\paragraph{Case 1:} $v$ is a 1-vertex.\\ 
Let $w_1$ denote the target of the unique downward edge emanating from
$v$ as in the following figure, where (here and in all subsequent
figures) each edge is labelled by the
probability of reaching it from $v$.
\begin{center}
\input {EPS/redge1.pstex_t}
\end{center}
In this case, $\EXP(v) = 1+\EXP(w_1)$. In the setup presented above,
we have
$\lambda_1=1$, $c=1$, $\Delta_1(w_1) \ge 1$, and $\Delta(w_1) \ge 1$,
thus (\ref{ind}) is implied by
\begin{equation} \label{eq1}
\alpha+\beta\ge 1.
\end{equation} 

\paragraph{Case 2:} $v$ is a 2-vertex.\\
Let $w_1$ and $w_2$ denote the targets of the two downward edges
emanating from $v$, where $w_2$ is lower than $w_1$.
\begin{center}
\input {EPS/redge2.pstex_t}
\end{center}
We have
\[
\EXP(v)\ =\ 1+\frac{1}{2}\EXP({w_1}) + \frac{1}{2}\EXP({w_2}),
\]
hence we need to require that
\[
\frac{\alpha}{2}\Delta_1(w_1) + 
\frac{\alpha}{2}\Delta_1(w_2) + 
\frac{\beta}{2}\Delta(w_1) + 
\frac{\beta}{2}\Delta(w_2)\ \ge\ 1.
\]
Note that $\Delta(w_2)>\Delta(w_1)\ge 1$.

\paragraph{Case 2.a:} $\Delta(w_2)\ge 4$ (as in the preceding figure).\\
Ignoring the effect of the $\Delta_1(w_j)$'s, it suffices to require
that
\[
\frac{\beta}{2}\Delta(w_1) + \frac{\beta}{2}\Delta(w_2) \ge 1,
\]
which will follow if
\begin{equation} \label{eq2}
\beta\ge \frac{2}{5}.
\end{equation}

\paragraph{Case 2.b.i:} $\Delta(w_2) = 3$ and one of the two vertices above
$w_2$ and below $v$ is a 1-vertex.\\
In this case $\Delta_1(w_2)\ge 1$ and $\Delta(w_1)+\Delta(w_2)\ge 4$,
so (\ref{ind}) is implied by
\begin{equation} \label{eq3}
\frac{1}{2}\alpha+2\beta\ge 1.
\end{equation}

\paragraph{Case 2.b.ii:} $\Delta(w_2) = 3$ and the two vertices between
$v$ and $w_2$ are 2-vertices.
Denote the second intermediate vertex as $v'$.
We may assume that $v'$ is reachable from $v$ (that is, from $w_1$),
otherwise we can ignore it and reduce the situation to Case 2.c
treated below
(by choosing another ordering of the vertices producing the same
oriented graph)%
. Three subcases can arise.

First, assume that none of the three edges that emanate from $w_1$ and
$v'$ further down reaches $w_2$. Denote by $x,y$ the two downward
neighbors of $v'$ and by $z$ the downward neighbor of $w_1$ other than
$v'$. The vertices $x,y,z$ need not be distinct (except that $x\ne y$),
but none of them coincides with $w_2$.
\begin{center}
\input EPS-A/redge2bii1.pstex_t
\end{center}
We have here $c=7/4$.

To make the analysis simpler to follow visually, we present it in a
table. Each row denotes one of the target vertices $w_2,x,y,z$,
`multiplied' by the probability of reaching it from $v$. The left 
(resp., right) column denotes a lower bound on the corresponding quantities 
$\Delta_1(\cdot)$ (resp., $\Delta(\cdot)$). 
To obtain an inequality that implies (2), one has to multiply
each entry in the left (resp., right) column by the row probability
times $\alpha$ (resp., times $\beta$), and require that the sum of all
these terms be $\ge c$.
\begin{center}
\begin{tabular}{|l|c|c|}
& $\alpha \Delta_1$ & $\beta \Delta$ \\
\hline
$1/2 w_2$ & 0 & 3 \\
$1/8 x$ & 0 & 4 \\
$1/8 y$ & 0 & 5 \\
$1/4 z$ & 0 & 4 \\
\hline
\end{tabular}
\end{center}
Note the following: (a) We do not assume that the rows represent 
distinct vertices (in fact, $x=z$ is implicit in the table); 
this does not cause any problem
in applying the rule for deriving an inequality from the table.
(b) We 
have
to squeeze the vertices so as to make the resulting
inequality as sharp (and difficult to satisfy) as possible; 
thus we made one of $x,y$ the farthest
vertex, because making $z$ the farthest vertex would have made the
inequality easier to satisfy.

We thus obtain
\[
\left(\frac{3}{2}+\frac{4}{8}+\frac{5}{8}+\frac{4}{4}\right)\beta\ge
\frac{7}{4} ,
\]
or
\begin{equation} \label{app:eq4}
\beta\ge \frac{14}{29} .
\end{equation}

Next, assume that $w_2$ is 
connected to
$v'$. In this case $w_2$ is a
1-vertex, and we extend the configuration to include its unique downward
neighbor $w_3$.
\begin{center}
\input EPS-A/redge2bii2.pstex_t
\end{center}
Let $x$ denote the other downward neighbor of $v'$ and 
let $y$ denote the other downward neighbor of $w_1$. 
In the following table, the `worst' case is to make $w_3$ and $y$ 
coincide, and make $x$ the farthest vertex.
\begin{center}
\begin{tabular}{|l|c|c|}
& $\alpha \Delta_1$ & $\beta \Delta$ \\
\hline
$5/8 w_3$ & 1 & 4 \\
$1/8 x$ & 1 & 5 \\
$1/4 y$ & 1 & 4 \\
\hline
\end{tabular}
\end{center}
We then obtain
\[
\alpha + \left(\frac{20}{8}+\frac{5}{8}+\frac{4}{4}\right)\beta\ge
\frac{19}{8},
\]
or
\begin{equation} \label{app:eq5}
\alpha + \frac{33}{8}\beta\ge \frac{19}{8}.
\end{equation}

Finally, assume that $w_2$ is 
connected to
$w_1$. Here too $w_2$ is a
1-vertex, and we extend the configuration to include its unique downward
neighbor $w_3$. 
\begin{center}
\input EPS-A/redge2bii3.pstex_t
\end{center}
Denoting by $x,y$ the two downward neighbors of $v'$, our 
table and resulting inequality become
\begin{equation} \label{app:eq6}
\begin{tabular}{|l|c|c|}
& $\alpha \Delta_1$ & $\beta \Delta$ \\
\hline
$3/4 w_3$ & 1 & 4 \\
$1/8 x$ & 1 & 4 \\
$1/8 y$ & 1 & 5 \\
\hline
\end{tabular}
\hspace{2cm}
\alpha + \frac{33}{8}\beta\ge \frac{5}{2} ,
\end{equation}
which, by the way, is stronger than (\ref{app:eq5}).

\paragraph{Case 2.c:} $\Delta(w_2) = 2$.
Hence, the only remaining case is that $w_1$ and $w_2$ are the
two vertices immediately following $v$.

\paragraph{Case 2.c.i:} $w_1$ is a 1-vertex (whose other upward neighbor
lies above $v$). Its unique downward edge ends at some vertex 
which is either $w_2$ or lies below $w_2$. 

Assume first that this vertex coincides with $w_2$, which makes $w_2$ a
1-vertex, whose unique downward neighbor is denoted as $v'$.
The local structure, table, and inequality are
\begin{equation} \label{app:eq7}
\input EPS-A/redge2ci1.pstex_t
%
\hspace{1.5cm}
\begin{tabular}{|l|c|c|}
& $\alpha \Delta_1$ & $\beta \Delta$ \\
\hline
$v'$ & 2 & 3 \\
\hline
\end{tabular}
\hspace{1.5cm}
2\alpha + 3\beta\ge \frac{5}{2} .
\end{equation}

Suppose next that the downward neighbor $w_3$ of $w_1$ lies below $w_2$.
We get
\begin{equation} \label{app:eq8}
\input EPS-A/redge2ci2.pstex_t
%
\hspace{1.5cm}
\begin{tabular}{|l|c|c|}
& $\alpha \Delta_1$ & $\beta \Delta$ \\
\hline
$1/2 w_2$ & 1 & 2 \\
$1/2 w_3$ & 1 & 3 \\
\hline
\end{tabular}
\hspace{1.5cm}
\alpha + \frac{5}{2}\beta\ge \frac{3}{2} .
\end{equation}

\paragraph{Case 2.c.ii:} $w_1$ is a 2-vertex, both of whose downward 
neighbors lie strictly below $w_2$. Denote these neighbors as $w_3,w_4$, 
with $w_3$ lying above $w_4$.
\begin{center}
\input EPS-A/redge2cii.pstex_t
\end{center}

We may assume that $\Delta(w_3)=3$ (i.\,e., there is no vertex between~$w_2$ 
and~$w_3$), since the case
$\Delta(w_3)\ge 4$ 
is already covered by~(\ref{app:eq4}).
%

\paragraph{Case 2.c.ii.1:} $w_2$ is a 1-vertex. 
Then 
the table and inequality become
\begin{equation} \label{app:eq9}
\begin{tabular}{|l|c|c|}
& $\alpha \Delta_1$ & $\beta \Delta$ \\
\hline
$1/2 w_2$ & 0 & 2 \\
$1/4 w_3$ & 1 & 3 \\
$1/4 w_4$ & 1 & 4 \\
\hline
\end{tabular}
\hspace{2cm}
\frac{1}{2}\alpha + \frac{11}{4}\beta\ge \frac{3}{2} .
\end{equation}

\paragraph{Case 2.c.ii.2:} $w_2$ is a 2-vertex but $w_3$ is a 1-vertex. 
Then $w_3$ (which satisfies $\Delta(w_3)=3$) is
connected either to $w_2$ or to a vertex above $v$. 
In the former
case, let $x$ denote the other downward neighbor of $w_2$, and  
let $y$ denote the unique downward neighbor of $w_3$. The local
structure looks like this (with $x,y,w_4$ not necessarily distinct,
but they all are below $w_3$ due to $\Delta(w_3)=3$): 
\begin{center}
\input EPS-A/redge2cii1.pstex_t
\end{center}
%
The (worst) table and inequality are
\begin{equation} \label{app:eq11}
\begin{tabular}{|l|c|c|}
& $\alpha \Delta_1$ & $\beta \Delta$ \\
\hline
$1/4 x$ & 1 & 4 \\
$1/2 y$ & 1 & 4 \\
$1/4 w_4$ & 1 & 5 \\
\hline
\end{tabular}
\hspace{2cm}
\alpha + \frac{17}{4}\beta\ge \frac{5}{2} .
\end{equation}

The next case is where the other upward neighbor of $w_3$ lies above
$v$. Let $x,y$ denote the two downward neighbors of $w_2$, and let $z$ 
denote the unique downward neighbor of $w_3$.
(Again, $x,y,z,w_4$ need not be distinct, but $x\ne y$ and they all are below $w_3$ due
  to $\Delta(w_3)=3$.) The local structure is:
\begin{center}
\input EPS-A/redge2cii2.pstex_t
\end{center}
%
%
The (worst) table and inequality become
\begin{equation} \label{app:eq14}
\begin{tabular}{|l|c|c|}
& $\alpha \Delta_1$ & $\beta \Delta$ \\
\hline
$1/4 x$ & 1 & 4 \\
$1/4 y$ & 1 & 5 \\
$1/4 z$ & 1 & 4 \\
$1/4 w_4$ & 1 & 5 \\
\hline
\end{tabular}
\hspace{2cm}
\alpha + \frac{9}{2}\beta\ge \frac{9}{4} .
\end{equation}

\paragraph{Case 2.c.ii.3:} Both $w_2$ and $w_3$ are 2-vertices. 
We have to consider the following type of configuration (where 
$x,y,z,t,w_4$ need not all be distinct, but $x\not= y$ and $z\not= t$, 
and we may assume $x\not=t$, $y\not= z$; also, because $\Delta(w_3)=3$,
both~$x$ and~$y$ are lower than~$w_3$):
\begin{center}
\input EPS-A/redge2cii3.pstex_t
\end{center}
Intuitively, a worst table is obtained by `squeezing' $x,y,z,t$, and 
$w_4$ as much to the left as possible, placing two of them at distance 4
from $v$, two at distance 5, and one at distance 6. However, squeezing
them this way will make some pairs of them coincide and form 1-vertices,
which will affect the resulting tables and inequalities.

Suppose first that among the three `heavier' targets $x,y,w_4$, at most
one lies at distance 4 from $v$. The worst table and the associated
inequality are (recall that $x\ne y$):
\begin{equation} \label{app:eq14a}
\begin{tabular}{|l|c|c|}
& $\alpha \Delta_1$ & $\beta \Delta$ \\
\hline
$1/4 x$ & 0 & 4 \\
$1/4 y$ & 0 & 5 \\
$1/8 z$ & 0 & 4 \\
$1/8 t$ & 0 & 6 \\
$1/4 w_4$ & 0 & 5 \\
\hline
\end{tabular}
\hspace{2cm}
\frac{19}{4}\beta\ge \frac{9}{4} .
\end{equation}

%

Suppose then that among $\{w_4,x,y\}$, two are at distance 4 from $v$,
say $w_4$ and $y$. Then $w_4=y$ is a 1-vertex, and we denote by $w$ its
unique downward neighbor. The local structure is:
\begin{center}
  \input EPS-A/redge2cii4.pstex_t
\end{center}
Two equally worst tables, and the resulting common inequality are

\begin{equation} \label{app:eq:new2}
\begin{tabular}{|l|c|c|}
& $\alpha \Delta_1$ & $\beta \Delta$ \\
\hline
$1/4 x$   & 1 & 5 \\
$1/8 z$   & 1 & 6 \\
$1/8 t$   & 1 & 7 \\
$1/2 w$ & 1 & 5 \\
\hline
\end{tabular}
\hspace{1.5cm}
\begin{tabular}{|l|c|c|}
& $\alpha \Delta_1$ & $\beta \Delta$ \\
\hline
$1/4 x$   & 1 & 6 \\
$1/8 z$   & 1 & 6 \\
$1/8 t$   & 1 & 5 \\
$1/2 w$ & 1 & 5 \\
\hline
\end{tabular}
\hspace{1.5cm}
  \alpha+\frac{43}{8}\beta\ge \frac{11}{4}\ .
\end{equation}

%
%

\paragraph{Case 2.c.iii:} $w_1$ is a 2-vertex that reaches $w_2$; that 
is, one of its downward neighbors, say $w_3$, coincides with $w_2$. Then
$w_2$ is a 1-vertex, and we denote by $x$ its unique downward neighbor.
\begin{center}
\input EPS-A/redge2ciii.pstex_t
\end{center}
A crucial observation is that $x$ cannot be equal to $w_4$. Indeed, if 
they were equal, then $w_4$ would be a 1-vertex. 
\begin{center}
\input EPS-A/redge2ciiix.pstex_t
\end{center}
In this case, cutting the edge graph $G$ of $P$
at the downward edge emanating from $w_4$ and at the edge entering $v$
would have disconnected $G$, contradicting the fact that $G$
is 3-connected%
.

We first dispose of the case where $x$ lies lower than $w_4$. 
The table and inequality are
\begin{equation} \label{app:eq15}
\begin{tabular}{|l|c|c|}
& $\alpha \Delta_1$ & $\beta \Delta$ \\
\hline
$3/4 x$ & 1 & 4 \\
$1/4 w_4$ & 1 & 3 \\
\hline
\end{tabular}
\hspace{2cm}
\alpha + \frac{15}{4}\beta\ge \frac{9}{4} .
\end{equation}
In what follows we thus assume that $x$ lies above $w_4$.

\paragraph{Case 2.c.iii.1:} $x$ is a 1-vertex that precedes $w_4$.
Suppose first that $w_4$ is the unique downward neighbor of $x$. Then
$w_4$ is a 1-vertex, and we denote its unique downward neighbor by $z$.
The local structure, table and inequality are: 
\begin{equation} \label{app:eq16}
\input EPS-A/redge2ciii1a.pstex_t
\hspace{0.5cm}
\begin{tabular}{|l|c|c|}
& $\alpha \Delta_1$ & $\beta \Delta$ \\
\hline
$z$ & 3 & 5 \\
\hline
\end{tabular}
\hspace{1cm}
3\alpha + 5\beta\ge 4 .
\end{equation}

Suppose next that the unique downward neighbor $y$ of $x$ is not $w_4$.
The local structure, table and inequality look like this 
($y$ is drawn above $w_4$ because this yields a sharper inequality):
\begin{equation} \label{app:eq17}
\input EPS-A/redge2ciii1b.pstex_t
\hspace{0.5cm}
\begin{tabular}{|l|c|c|}
& $\alpha \Delta_1$ & $\beta \Delta$ \\
\hline
$3/4 y$ & 2 & 4 \\
$1/4 w_4$ & 2 & 5 \\
\hline
\end{tabular}
\hspace{1cm}
2\alpha + \frac{17}{4}\beta\ge 3 .
\end{equation}

\paragraph{Case 2.c.iii.2:} $x$ is a 2-vertex that precedes $w_4$.
This subcase splits into several subcases, where we assume,
respectively, that $\Delta(w_4)\ge 6$,
$\Delta(w_4) = 4$, and $\Delta(w_4) = 5$.

\paragraph{Case 2.c.iii.2(a).} 
Suppose first that $\Delta(w_4)\ge 6$. The configuration
looks like this:
\begin{center}
\input EPS-A/redge2ciii2.pstex_t
\end{center}
The table and inequality are
\begin{equation} \label{app:eq18}
\begin{tabular}{|l|c|c|}
& $\alpha \Delta_1$ & $\beta \Delta$ \\
\hline
$3/4 x$ & 1 & 3 \\
$1/4 w_4$ & 1 & 6 \\
\hline
\end{tabular}
\hspace{2cm}
\alpha + \frac{15}{4}\beta\ge \frac{9}{4} .
\end{equation}
Note that this is the same inequality as (\ref{app:eq15}).

\paragraph{Case 2.c.iii.2(b).} 
Suppose next that $\Delta(w_4)=4$, and that one of the downward neighbors 
of $x$ is $w_4$.
Let $z$ denote the other downward neighbor. $w_4$ is a 1-vertex, and we
denote by $w$ its unique downward neighbor.
\begin{center}
\input EPS-A/redge2ciii2a.pstex_t
\end{center}
The 3-connectivity of the edge graph of $P$ implies, as above, that
$w\ne z$. 
Since we assume that $\Delta(w_4)=4$, $z$ 
also lies below $w_4$, and the table and inequality are
\begin{equation} \label{app:eq19}
\begin{tabular}{|l|c|c|}
& $\alpha \Delta_1$ & $\beta \Delta$ \\
\hline
$5/8 w$ & 2 & 5 \\
$3/8 z$ & 2 & 6 \\
\hline
\end{tabular}
\hspace{2cm}
2\alpha + \frac{43}{8}\beta\ge \frac{29}{8} .
\end{equation}

Suppose next that $\Delta(w_4)=4$ and $w_4$ is not a downward 
neighbor of $x$. Denote those two neighbors as $w$ and $z$, both of
which lie lower than $w_4$, by assumption, and are clearly
distinct. The configuration, table and inequality look like this:
\begin{equation} \label{app:eq20}
\input EPS-A/redge2ciii2b.pstex_t
\hspace{0.5cm}
\begin{tabular}{|l|c|c|}
& $\alpha \Delta_1$ & $\beta \Delta$ \\
\hline
$1/4 w_4$ & 1 & 4 \\
$3/8 w$ & 1 & 5 \\
$3/8 z$ & 1 & 6 \\
\hline
\end{tabular}
\hspace{1cm}
\alpha + \frac{41}{8}\beta\ge 3 .
\end{equation}

\paragraph{Case 2.c.iii.2(c).} 
It remains to consider the case $\Delta(w_4)=5$. Let $z$ denote the 
unique vertex lying between $x$ and $w_4$. We may assume that $z$ is 
connected to $x$, for otherwise $z$ is not reachable from $v$, and 
we might as well reduce this case to the case $\Delta(w_4)=4$ just treated.

Consider first the subcase where the other downward neighbor of
$x$ is $w_4$ itself. Then $w_4$ is a 1-vertex, and we denote by
$w$ its unique downward neighbor. This subcase splits further into 
two subcases: First, assume that $z$ is a 1-vertex, and let $y$
denote its unique downward neighbor. Clearly, $y$ must lie below
$w_4$ (it may coincide with or precede $w$). The configuration
looks like this:
\begin{center}
\input EPS-A/redge2civ1a.pstex_t
\end{center}
The table and inequality are
\begin{equation} \label{app:eq21}
\begin{tabular}{|l|c|c|}
& $\alpha \Delta_1$ & $\beta \Delta$ \\
\hline
$3/8 y$ & 3 & 6 \\
$5/8 w$ & 3 & 6 \\
\hline
\end{tabular}
\hspace{2cm}
3\alpha + 6\beta\ge 4 .
\end{equation}

In the other subcase, $z$ is a 2-vertex; we denote its two
downward neighbors as $y$ and $t$. The vertices $w,y,t$ all lie
below $w_4$ and may appear there in any order (except that $w\ne t$). 
The configuration looks like this:
\begin{center}
\input EPS-A/redge2civ1b.pstex_t
\end{center}
The table and inequality are 
\begin{equation} \label{app:eq22}
\begin{tabular}{|l|c|c|}
& $\alpha \Delta_1$ & $\beta \Delta$ \\
\hline
$3/16 y$ & 2 & 6 \\
$3/16 t$ & 2 & 7 \\
$5/8 w$ & 2 & 6 \\
\hline
\end{tabular}
\hspace{2cm}
2\alpha + \frac{99}{16}\beta\ge 4 .
\end{equation}

Consider next the subcase where $w_4$ is not a downward neighbor of
$x$. Denote the other downward neighbor of $x$ as $y$, which lies
strictly below $w_4$.
This subcase splits into three subcases. First, assume that
$z$ is a 1-vertex, and denote its unique downward neighbor as $w$.
The configuration looks like this:
\begin{center}
\input EPS-A/redge2civ2a.pstex_t
\end{center}
The table and inequality are
\begin{equation} \label{app:eq23}
\begin{tabular}{|l|c|c|}
& $\alpha \Delta_1$ & $\beta \Delta$ \\
\hline
$1/4 w_4$ & 2 & 5 \\
$3/8 y$ & 2 & 6 \\
$3/8 w$ & 2 & 5 \\
\hline
\end{tabular}
\hspace{2cm}
2\alpha + \frac{43}{8}\beta\ge \frac{27}{8} .
\end{equation}

Second, assume that $z$ is a 2-vertex, so that none of its two
downward neighbors is $w_4$. Denote these neighbors as $w$ and
$t$. All three vertices $y,t,w$ lie strictly below $w_4$, and
$w\ne t$. The configuration looks like this:
\begin{center}
\input EPS-A/redge2civ2b.pstex_t
\end{center}
The table and inequality are
\begin{equation} \label{app:eq24}
\begin{tabular}{|l|c|c|}
& $\alpha \Delta_1$ & $\beta \Delta$ \\
\hline
$1/4 w_4$ & 1 & 5 \\
$3/8 y$ & 1 & 6 \\
$3/16 w$ & 1 & 6 \\
$3/16 t$ & 1 & 7 \\
\hline
\end{tabular}
\hspace{2cm}
\alpha + \frac{95}{16}\beta\ge \frac{27}{8} .
\end{equation}

Finally, assume that $z$ is a 2-vertex, so that one of its two
downward neighbors is $w_4$. Denote the other neighbor as $w$. In
this case $w_4$ is a 1-vertex, and we denote its unique downward
neighbor as $t$. All three vertices $y,t,w$ lie strictly below $w_4$.
The configuration looks like this:
\begin{center}
\input EPS-A/redge2civ2c.pstex_t
\end{center}
The table and inequality are
\begin{equation} \label{app:eq25}
\begin{tabular}{|l|c|c|}
& $\alpha \Delta_1$ & $\beta \Delta$ \\
\hline
$3/8 y$ & 2 & 6 \\
$3/16 w$ & 2 & 7 \\
$7/16 t$ & 2 & 6 \\
\hline
\end{tabular}
\hspace{2cm}
2\alpha + \frac{99}{16}\beta\ge \frac{61}{16} ,
\end{equation}
which, by the way, is weaker than (\ref{app:eq22}).

\medskip

This completes the case distinction. Thus~\eqref{eq:generic} holds for every pair $(\alpha,\beta)$ that satisfies~\eqref{eq1}--\eqref{app:eq25}. In particular, it holds for the pair $(\alpha,\beta)=(\frac{46}{87},\frac{42}{87})$, which (as discussed at the beginning of the proof) yields the upper bound $\frac{130}{87}<1.4943$ on the linearity coefficient of random edge.

\end{proof}

\medskip
\noindent{\bf Discussion.}
(1) The analysis has used (twice) the fact that $G$ is 
a 3-connected graph. Without this assumption, the linearity 
coefficient becomes 13/8: A lower bound construction can be derived
from the figure shown in Case 2.c.iii, and an upper bound can be 
obtained along the same lines of the preceding proof, using 
a much shorter case analysis. 
It is interesting that the proof did not use at 
all the planarity of the polytope graph $G$.

(2) In an earlier phase of our work, we obtained the upper
bound of $3/2$ on the linearity coefficient, using a similar
but considerably shorter case analysis. 
Unfortunately, the lengthier case distinction  presented in the proof above is not just a refinement of that shorter one (which is the reason for presenting only the lengthier proof). 
The proof indicates that the problem probably is  far from admitting a 
clean and simple solution -- at least using this approach. 
Of course, it would be interesting to find an alternative 
simpler way of attacking the problem.

(3) The solution
$(\alpha,\beta)=(\frac{46}{87},\frac{42}{87})$ satisfies 
(\ref{app:eq4}) and (\ref{app:eq16}) with equality. If we examine 
the configuration corresponding to (\ref{app:eq4}) and expand 
it further, we can replace (\ref{app:eq4}) by better inequalities,
which result in a (slightly) improved bound on the linearity 
coefficient, at the cost of lengthening further our case analysis.
This refinement process can continue for a few more steps, as we have verified.
We have no idea whether this iterative refinement process ever
converges to some critical configuration, whose further expansion
does not improve the bound, and which is then likely to yield
a tight bound on the linearity coefficient.


%% file: EPS/backbone.pstex_t
\begin{picture}(0,0)%
\includegraphics{backbone.pstex}%
\end{picture}%
\setlength{\unitlength}{4144sp}%
\begingroup\makeatletter\ifx\SetFigFont\undefined
\def\x#1#2#3#4#5#6#7\relax{\def\x{#1#2#3#4#5#6}}%
\expandafter\x\fmtname xxxxxx\relax \def\y{splain}%
\ifx\x\y   
\gdef\SetFigFont#1#2#3{%
  \ifnum #1<17\tiny\else \ifnum #1<20\small\else
  \ifnum #1<24\normalsize\else \ifnum #1<29\large\else
  \ifnum #1<34\Large\else \ifnum #1<41\LARGE\else
     \huge\fi\fi\fi\fi\fi\fi
  \csname #3\endcsname}%
\else
\gdef\SetFigFont#1#2#3{\begingroup
  \count@#1\relax \ifnum 25<\count@\count@25\fi
  \def\x{\endgroup\@setsize\SetFigFont{#2pt}}%
  \expandafter\x
    \csname \romannumeral\the\count@ pt\expandafter\endcsname
    \csname @\romannumeral\the\count@ pt\endcsname
  \csname #3\endcsname}%
\fi
\fi\endgroup
\begin{picture}(3600,936)(631,-283)
\put(2521,254){\makebox(0,0)[b]{\smash{\SetFigFont{8}{9.6}{rm}{\color[rgb]{0,0,0}$\cdots$}%
}}}
\put(1261,-106){\makebox(0,0)[b]{\smash{\SetFigFont{8}{9.6}{rm}{\color[rgb]{0,0,0}$v_{2k-7}$}%
}}}
\put(1891,-106){\makebox(0,0)[b]{\smash{\SetFigFont{8}{9.6}{rm}{\color[rgb]{0,0,0}$v_{2k-9}$}%
}}}
\put(2431,-106){\makebox(0,0)[b]{\smash{\SetFigFont{8}{9.6}{rm}{\color[rgb]{0,0,0}$v_{2k-11}$}%
}}}
\put(2071,569){\makebox(0,0)[b]{\smash{\SetFigFont{8}{9.6}{rm}{\color[rgb]{0,0,0}$v_{2k-10}$}%
}}}
\put(1531,569){\makebox(0,0)[b]{\smash{\SetFigFont{8}{9.6}{rm}{\color[rgb]{0,0,0}$v_{2k-8}$}%
}}}
\put(991,569){\makebox(0,0)[b]{\smash{\SetFigFont{8}{9.6}{rm}{\color[rgb]{0,0,0}$v_{2k-6}$}%
}}}
\put(631,-16){\makebox(0,0)[rb]{\smash{\SetFigFont{8}{9.6}{rm}{\color[rgb]{0,0,0}$v_{2k-5} = v_{\max}$}%
}}}
\put(4231,-16){\makebox(0,0)[lb]{\smash{\SetFigFont{8}{9.6}{rm}{\color[rgb]{0,0,0}$v_0=v_{\min}$}%
}}}
\put(3061,-106){\makebox(0,0)[b]{\smash{\SetFigFont{8}{9.6}{rm}{\color[rgb]{0,0,0}$v_3$}%
}}}
\put(2701,569){\makebox(0,0)[b]{\smash{\SetFigFont{8}{9.6}{rm}{\color[rgb]{0,0,0}$v_4$}%
}}}
\put(3241,569){\makebox(0,0)[b]{\smash{\SetFigFont{8}{9.6}{rm}{\color[rgb]{0,0,0}$v_2$}%
}}}
\put(3601,-106){\makebox(0,0)[b]{\smash{\SetFigFont{8}{9.6}{rm}{\color[rgb]{0,0,0}$v_1$}%
}}}
\end{picture}

%% file: EPS/re.pstex_t
\begin{picture}(0,0)%
\includegraphics{re.pstex}%
\end{picture}%
\setlength{\unitlength}{4144sp}%
\begingroup\makeatletter\ifx\SetFigFont\undefined
\def\x#1#2#3#4#5#6#7\relax{\def\x{#1#2#3#4#5#6}}%
\expandafter\x\fmtname xxxxxx\relax \def\y{splain}%
\ifx\x\y   
\gdef\SetFigFont#1#2#3{%
  \ifnum #1<17\tiny\else \ifnum #1<20\small\else
  \ifnum #1<24\normalsize\else \ifnum #1<29\large\else
  \ifnum #1<34\Large\else \ifnum #1<41\LARGE\else
     \huge\fi\fi\fi\fi\fi\fi
  \csname #3\endcsname}%
\else
\gdef\SetFigFont#1#2#3{\begingroup
  \count@#1\relax \ifnum 25<\count@\count@25\fi
  \def\x{\endgroup\@setsize\SetFigFont{#2pt}}%
  \expandafter\x
    \csname \romannumeral\the\count@ pt\expandafter\endcsname
    \csname @\romannumeral\the\count@ pt\endcsname
  \csname #3\endcsname}%
\fi
\fi\endgroup
\begin{picture}(5519,1299)(-11,-497)
\put(541,-241){\makebox(0,0)[b]{\smash{\SetFigFont{8}{9.6}{rm}{\color[rgb]{0,0,0}$64$}%
}}}
\put(901,-241){\makebox(0,0)[b]{\smash{\SetFigFont{8}{9.6}{rm}{\color[rgb]{0,0,0}$32$}%
}}}
\put(1351,-241){\makebox(0,0)[b]{\smash{\SetFigFont{8}{9.6}{rm}{\color[rgb]{0,0,0}$16$}%
}}}
\put(2836,-466){\makebox(0,0)[b]{\smash{\SetFigFont{8}{9.6}{rm}{\color[rgb]{0,0,0}$8$}%
}}}
\put(541,164){\makebox(0,0)[b]{\smash{\SetFigFont{8}{9.6}{rm}{\color[rgb]{0,0,0}$64$}%
}}}
\put(991,164){\makebox(0,0)[b]{\smash{\SetFigFont{8}{9.6}{rm}{\color[rgb]{0,0,0}$96$}%
}}}
\put(1306,164){\makebox(0,0)[b]{\smash{\SetFigFont{8}{9.6}{rm}{\color[rgb]{0,0,0}$112$}%
}}}
\put(1576,164){\makebox(0,0)[b]{\smash{\SetFigFont{8}{9.6}{rm}{\color[rgb]{0,0,0}$56$}%
}}}
\put(1846,164){\makebox(0,0)[b]{\smash{\SetFigFont{8}{9.6}{rm}{\color[rgb]{0,0,0}$64$}%
}}}
\put(1711,434){\makebox(0,0)[b]{\smash{\SetFigFont{8}{9.6}{rm}{\color[rgb]{0,0,0}$56$}%
}}}
\put(2161,164){\makebox(0,0)[b]{\smash{\SetFigFont{8}{9.6}{rm}{\color[rgb]{0,0,0}$120$}%
}}}
\put(2566,164){\makebox(0,0)[b]{\smash{\SetFigFont{8}{9.6}{rm}{\color[rgb]{0,0,0}$60$}%
}}}
\put(2926,164){\makebox(0,0)[b]{\smash{\SetFigFont{8}{9.6}{rm}{\color[rgb]{0,0,0}$90$}%
}}}
\put(3196,164){\makebox(0,0)[b]{\smash{\SetFigFont{8}{9.6}{rm}{\color[rgb]{0,0,0}$90$}%
}}}
\put(3511,164){\makebox(0,0)[b]{\smash{\SetFigFont{8}{9.6}{rm}{\color[rgb]{0,0,0}$120$}%
}}}
\put(3871,389){\makebox(0,0)[b]{\smash{\SetFigFont{8}{9.6}{rm}{\color[rgb]{0,0,0}$60$}%
}}}
\put(3916, 29){\makebox(0,0)[b]{\smash{\SetFigFont{8}{9.6}{rm}{\color[rgb]{0,0,0}$60$}%
}}}
\put(4276,164){\makebox(0,0)[b]{\smash{\SetFigFont{8}{9.6}{rm}{\color[rgb]{0,0,0}$90$}%
}}}
\put(4546,164){\makebox(0,0)[b]{\smash{\SetFigFont{8}{9.6}{rm}{\color[rgb]{0,0,0}$98$}%
}}}
\put(4951,-286){\makebox(0,0)[b]{\smash{\SetFigFont{8}{9.6}{rm}{\color[rgb]{0,0,0}$49$}%
}}}
\put(5086,164){\makebox(0,0)[b]{\smash{\SetFigFont{8}{9.6}{rm}{\color[rgb]{0,0,0}$79$}%
}}}
\put(5356,164){\makebox(0,0)[b]{\smash{\SetFigFont{8}{9.6}{rm}{\color[rgb]{0,0,0}$128$}%
}}}
\put(4861,-16){\makebox(0,0)[b]{\smash{\SetFigFont{8}{9.6}{rm}{\color[rgb]{0,0,0}$49$}%
}}}
\put(811,-16){\makebox(0,0)[lb]{\smash{\SetFigFont{8}{9.6}{rm}{\color[rgb]{0,0,0}$32$}%
}}}
\put(1711,-16){\makebox(0,0)[lb]{\smash{\SetFigFont{8}{9.6}{rm}{\color[rgb]{0,0,0}$8$}%
}}}
\put(2836,-16){\makebox(0,0)[lb]{\smash{\SetFigFont{8}{9.6}{rm}{\color[rgb]{0,0,0}$30$}%
}}}
\put(3241,-151){\makebox(0,0)[lb]{\smash{\SetFigFont{8}{9.6}{rm}{\color[rgb]{0,0,0}$30$}%
}}}
\put(4096,164){\makebox(0,0)[rb]{\smash{\SetFigFont{8}{9.6}{rm}{\color[rgb]{0,0,0}$30$}%
}}}
\put(1171,-16){\makebox(0,0)[lb]{\smash{\SetFigFont{8}{9.6}{rm}{\color[rgb]{0,0,0}$16$}%
}}}
\put(4411,434){\makebox(0,0)[b]{\smash{\SetFigFont{8}{9.6}{rm}{\color[rgb]{0,0,0}$30$}%
}}}
\put(2506,-226){\makebox(0,0)[b]{\smash{\SetFigFont{8}{9.6}{rm}{\color[rgb]{0,0,0}$60$}%
}}}
\end{picture}

%% file: EPS/redge1.pstex_t
\begin{picture}(0,0)%
\includegraphics{redge1.pstex}%
\end{picture}%
\setlength{\unitlength}{3947sp}%
\begingroup\makeatletter\ifx\SetFigFont\undefined
\def\x#1#2#3#4#5#6#7\relax{\def\x{#1#2#3#4#5#6}}%
\expandafter\x\fmtname xxxxxx\relax \def\y{splain}%
\ifx\x\y   
\gdef\SetFigFont#1#2#3{%
  \ifnum #1<17\tiny\else \ifnum #1<20\small\else
  \ifnum #1<24\normalsize\else \ifnum #1<29\large\else
  \ifnum #1<34\Large\else \ifnum #1<41\LARGE\else
     \huge\fi\fi\fi\fi\fi\fi
  \csname #3\endcsname}%
\else
\gdef\SetFigFont#1#2#3{\begingroup
  \count@#1\relax \ifnum 25<\count@\count@25\fi
  \def\x{\endgroup\@setsize\SetFigFont{#2pt}}%
  \expandafter\x
    \csname \romannumeral\the\count@ pt\expandafter\endcsname
    \csname @\romannumeral\the\count@ pt\endcsname
  \csname #3\endcsname}%
\fi
\fi\endgroup
\begin{picture}(1050,298)(526,-278)
\put(526,-136){\makebox(0,0)[lb]{\smash{\SetFigFont{12}{14.4}{rm}{\color[rgb]{0,0,0}$v$}%
}}}
\put(1576,-136){\makebox(0,0)[lb]{\smash{\SetFigFont{12}{14.4}{rm}{\color[rgb]{0,0,0}$w_1$}%
}}}
\put(976,-156){\makebox(0,0)[lb]{\smash{\SetFigFont{12}{14.4}{rm}{\color[rgb]{0,0,0}$1$}%
}}}
\end{picture}

%% file: EPS/redge2.pstex_t
\begin{picture}(0,0)%
\includegraphics{redge2.pstex}%
\end{picture}%
\setlength{\unitlength}{3158sp}%
\begingroup\makeatletter\ifx\SetFigFont\undefined
\def\x#1#2#3#4#5#6#7\relax{\def\x{#1#2#3#4#5#6}}%
\expandafter\x\fmtname xxxxxx\relax \def\y{splain}%
\ifx\x\y   
\gdef\SetFigFont#1#2#3{%
  \ifnum #1<17\tiny\else \ifnum #1<20\small\else
  \ifnum #1<24\normalsize\else \ifnum #1<29\large\else
  \ifnum #1<34\Large\else \ifnum #1<41\LARGE\else
     \huge\fi\fi\fi\fi\fi\fi
  \csname #3\endcsname}%
\else
\gdef\SetFigFont#1#2#3{\begingroup
  \count@#1\relax \ifnum 25<\count@\count@25\fi
  \def\x{\endgroup\@setsize\SetFigFont{#2pt}}%
  \expandafter\x
    \csname \romannumeral\the\count@ pt\expandafter\endcsname
    \csname @\romannumeral\the\count@ pt\endcsname
  \csname #3\endcsname}%
\fi
\fi\endgroup
\begin{picture}(3360,839)(526,-440)
\put(526,-136){\makebox(0,0)[lb]{\smash{\SetFigFont{10}{12.0}{rm}{\color[rgb]{0,0,0}$v$}%
}}}
\put(1576,-136){\makebox(0,0)[lb]{\smash{\SetFigFont{10}{12.0}{rm}{\color[rgb]{0,0,0}$w_1$}%
}}}
\put(3886,-158){\makebox(0,0)[lb]{\smash{\SetFigFont{10}{12.0}{rm}{\color[rgb]{0,0,0}$w_2$}%
}}}
\put(2221,224){\makebox(0,0)[lb]{\smash{\SetFigFont{10}{12.0}{rm}{\color[rgb]{0,0,0}$1/2$}%
}}}
\put(923,-391){\makebox(0,0)[lb]{\smash{\SetFigFont{10}{12.0}{rm}{\color[rgb]{0,0,0}$1/2$}%
}}}
\end{picture}

%% file: EPS-A/redge2bii1.pstex_t
\begin{picture}(0,0)%
\includegraphics{redge2bii1.pstex}%
\end{picture}%
\setlength{\unitlength}{3158sp}%
\begingroup\makeatletter\ifx\SetFigFont\undefined
\def\x#1#2#3#4#5#6#7\relax{\def\x{#1#2#3#4#5#6}}%
\expandafter\x\fmtname xxxxxx\relax \def\y{splain}%
\ifx\x\y   
\gdef\SetFigFont#1#2#3{%
  \ifnum #1<17\tiny\else \ifnum #1<20\small\else
  \ifnum #1<24\normalsize\else \ifnum #1<29\large\else
  \ifnum #1<34\Large\else \ifnum #1<41\LARGE\else
     \huge\fi\fi\fi\fi\fi\fi
  \csname #3\endcsname}%
\else
\gdef\SetFigFont#1#2#3{\begingroup
  \count@#1\relax \ifnum 25<\count@\count@25\fi
  \def\x{\endgroup\@setsize\SetFigFont{#2pt}}%
  \expandafter\x
    \csname \romannumeral\the\count@ pt\expandafter\endcsname
    \csname @\romannumeral\the\count@ pt\endcsname
  \csname #3\endcsname}%
\fi
\fi\endgroup
\begin{picture}(3300,1793)(466,-1324)
\put(526,-136){\makebox(0,0)[lb]{\smash{\SetFigFont{10}{12.0}{rm}{\color[rgb]{0,0,0}$v$}%
}}}
\put(3758,-406){\makebox(0,0)[lb]{\smash{\SetFigFont{10}{12.0}{rm}{\color[rgb]{0,0,0}$x$}%
}}}
\put(3758,-811){\makebox(0,0)[lb]{\smash{\SetFigFont{10}{12.0}{rm}{\color[rgb]{0,0,0}$y$}%
}}}
\put(3766,-1200){\makebox(0,0)[lb]{\smash{\SetFigFont{10}{12.0}{rm}{\color[rgb]{0,0,0}$z$}%
}}}
\put(1501,337){\makebox(0,0)[lb]{\smash{\SetFigFont{10}{12.0}{rm}{\color[rgb]{0,0,0}$1/2$}%
}}}
\put(908,-406){\makebox(0,0)[lb]{\smash{\SetFigFont{10}{12.0}{rm}{\color[rgb]{0,0,0}$1/2$}%
}}}
\put(1756,-398){\makebox(0,0)[lb]{\smash{\SetFigFont{10}{12.0}{rm}{\color[rgb]{0,0,0}$1/4$}%
}}}
\put(2236,-1275){\makebox(0,0)[lb]{\smash{\SetFigFont{10}{12.0}{rm}{\color[rgb]{0,0,0}$1/4$}%
}}}
\put(2873,-608){\makebox(0,0)[lb]{\smash{\SetFigFont{10}{12.0}{rm}{\color[rgb]{0,0,0}$1/8$}%
}}}
\put(2416,-788){\makebox(0,0)[lb]{\smash{\SetFigFont{10}{12.0}{rm}{\color[rgb]{0,0,0}$1/8$}%
}}}
\put(466,-1163){\makebox(0,0)[lb]{\smash{\SetFigFont{10}{12.0}{rm}{\color[rgb]{0,0,0}$c=7/4$}%
}}}
\put(2363,-181){\makebox(0,0)[lb]{\smash{\SetFigFont{10}{12.0}{rm}{\color[rgb]{0,0,0}$v'$}%
}}}
\put(1576,-136){\makebox(0,0)[lb]{\smash{\SetFigFont{10}{12.0}{rm}{\color[rgb]{0,0,0}$w_1$}%
}}}
\put(3031,-75){\makebox(0,0)[lb]{\smash{\SetFigFont{10}{12.0}{rm}{\color[rgb]{0,0,0}$w_2$}%
}}}
\end{picture}

%% file: EPS-A/redge2bii2.pstex_t
\begin{picture}(0,0)%
\includegraphics{redge2bii2.pstex}%
\end{picture}%
\setlength{\unitlength}{3158sp}%
\begingroup\makeatletter\ifx\SetFigFont\undefined
\def\x#1#2#3#4#5#6#7\relax{\def\x{#1#2#3#4#5#6}}%
\expandafter\x\fmtname xxxxxx\relax \def\y{splain}%
\ifx\x\y   
\gdef\SetFigFont#1#2#3{%
  \ifnum #1<17\tiny\else \ifnum #1<20\small\else
  \ifnum #1<24\normalsize\else \ifnum #1<29\large\else
  \ifnum #1<34\Large\else \ifnum #1<41\LARGE\else
     \huge\fi\fi\fi\fi\fi\fi
  \csname #3\endcsname}%
\else
\gdef\SetFigFont#1#2#3{\begingroup
  \count@#1\relax \ifnum 25<\count@\count@25\fi
  \def\x{\endgroup\@setsize\SetFigFont{#2pt}}%
  \expandafter\x
    \csname \romannumeral\the\count@ pt\expandafter\endcsname
    \csname @\romannumeral\the\count@ pt\endcsname
  \csname #3\endcsname}%
\fi
\fi\endgroup
\begin{picture}(3435,1793)(466,-1324)
\put(526,-136){\makebox(0,0)[lb]{\smash{\SetFigFont{10}{12.0}{rm}{\color[rgb]{0,0,0}$v$}%
}}}
\put(1501,337){\makebox(0,0)[lb]{\smash{\SetFigFont{10}{12.0}{rm}{\color[rgb]{0,0,0}$1/2$}%
}}}
\put(908,-406){\makebox(0,0)[lb]{\smash{\SetFigFont{10}{12.0}{rm}{\color[rgb]{0,0,0}$1/2$}%
}}}
\put(1756,-398){\makebox(0,0)[lb]{\smash{\SetFigFont{10}{12.0}{rm}{\color[rgb]{0,0,0}$1/4$}%
}}}
\put(2236,-1275){\makebox(0,0)[lb]{\smash{\SetFigFont{10}{12.0}{rm}{\color[rgb]{0,0,0}$1/4$}%
}}}
\put(2416,-788){\makebox(0,0)[lb]{\smash{\SetFigFont{10}{12.0}{rm}{\color[rgb]{0,0,0}$1/8$}%
}}}
\put(2491,-383){\makebox(0,0)[lb]{\smash{\SetFigFont{10}{12.0}{rm}{\color[rgb]{0,0,0}$1/8$}%
}}}
\put(3151,-391){\makebox(0,0)[lb]{\smash{\SetFigFont{10}{12.0}{rm}{\color[rgb]{0,0,0}$5/8$}%
}}}
\put(3713,-1253){\makebox(0,0)[lb]{\smash{\SetFigFont{10}{12.0}{rm}{\color[rgb]{0,0,0}$y$}%
}}}
\put(3720,-825){\makebox(0,0)[lb]{\smash{\SetFigFont{10}{12.0}{rm}{\color[rgb]{0,0,0}$x$}%
}}}
\put(466,-1163){\makebox(0,0)[lb]{\smash{\SetFigFont{10}{12.0}{rm}{\color[rgb]{0,0,0}$c=19/8$}%
}}}
\put(2341,-144){\makebox(0,0)[lb]{\smash{\SetFigFont{10}{12.0}{rm}{\color[rgb]{0,0,0}$v'$}%
}}}
\put(1576,-136){\makebox(0,0)[lb]{\smash{\SetFigFont{10}{12.0}{rm}{\color[rgb]{0,0,0}$w_1$}%
}}}
\put(3031,-75){\makebox(0,0)[lb]{\smash{\SetFigFont{10}{12.0}{rm}{\color[rgb]{0,0,0}$w_2$}%
}}}
\put(3901,-278){\makebox(0,0)[lb]{\smash{\SetFigFont{10}{12.0}{rm}{\color[rgb]{0,0,0}$w_3$}%
}}}
\end{picture}

%% file: EPS-A/redge2bii3.pstex_t
\begin{picture}(0,0)%
\includegraphics{redge2bii3.pstex}%
\end{picture}%
\setlength{\unitlength}{3158sp}%
\begingroup\makeatletter\ifx\SetFigFont\undefined
\def\x#1#2#3#4#5#6#7\relax{\def\x{#1#2#3#4#5#6}}%
\expandafter\x\fmtname xxxxxx\relax \def\y{splain}%
\ifx\x\y   
\gdef\SetFigFont#1#2#3{%
  \ifnum #1<17\tiny\else \ifnum #1<20\small\else
  \ifnum #1<24\normalsize\else \ifnum #1<29\large\else
  \ifnum #1<34\Large\else \ifnum #1<41\LARGE\else
     \huge\fi\fi\fi\fi\fi\fi
  \csname #3\endcsname}%
\else
\gdef\SetFigFont#1#2#3{\begingroup
  \count@#1\relax \ifnum 25<\count@\count@25\fi
  \def\x{\endgroup\@setsize\SetFigFont{#2pt}}%
  \expandafter\x
    \csname \romannumeral\the\count@ pt\expandafter\endcsname
    \csname @\romannumeral\the\count@ pt\endcsname
  \csname #3\endcsname}%
\fi
\fi\endgroup
\begin{picture}(3435,2005)(466,-1311)
\put(526,-136){\makebox(0,0)[lb]{\smash{\SetFigFont{10}{12.0}{rm}{\color[rgb]{0,0,0}$v$}%
}}}
\put(908,-406){\makebox(0,0)[lb]{\smash{\SetFigFont{10}{12.0}{rm}{\color[rgb]{0,0,0}$1/2$}%
}}}
\put(1756,-398){\makebox(0,0)[lb]{\smash{\SetFigFont{10}{12.0}{rm}{\color[rgb]{0,0,0}$1/4$}%
}}}
\put(3713,-1253){\makebox(0,0)[lb]{\smash{\SetFigFont{10}{12.0}{rm}{\color[rgb]{0,0,0}$y$}%
}}}
\put(3720,-825){\makebox(0,0)[lb]{\smash{\SetFigFont{10}{12.0}{rm}{\color[rgb]{0,0,0}$x$}%
}}}
\put(1606,562){\makebox(0,0)[lb]{\smash{\SetFigFont{10}{12.0}{rm}{\color[rgb]{0,0,0}$1/2$}%
}}}
\put(2768,-563){\makebox(0,0)[lb]{\smash{\SetFigFont{10}{12.0}{rm}{\color[rgb]{0,0,0}$1/8$}%
}}}
\put(3151,-391){\makebox(0,0)[lb]{\smash{\SetFigFont{10}{12.0}{rm}{\color[rgb]{0,0,0}$3/4$}%
}}}
\put(1906,142){\makebox(0,0)[lb]{\smash{\SetFigFont{10}{12.0}{rm}{\color[rgb]{0,0,0}$1/4$}%
}}}
\put(2610,-1096){\makebox(0,0)[lb]{\smash{\SetFigFont{10}{12.0}{rm}{\color[rgb]{0,0,0}$1/8$}%
}}}
\put(466,-1163){\makebox(0,0)[lb]{\smash{\SetFigFont{10}{12.0}{rm}{\color[rgb]{0,0,0}$c=5/2$}%
}}}
\put(2348,-181){\makebox(0,0)[lb]{\smash{\SetFigFont{10}{12.0}{rm}{\color[rgb]{0,0,0}$v'$}%
}}}
\put(1254,-83){\makebox(0,0)[lb]{\smash{\SetFigFont{10}{12.0}{rm}{\color[rgb]{0,0,0}$w_1$}%
}}}
\put(3031,-75){\makebox(0,0)[lb]{\smash{\SetFigFont{10}{12.0}{rm}{\color[rgb]{0,0,0}$w_2$}%
}}}
\put(3901,-278){\makebox(0,0)[lb]{\smash{\SetFigFont{10}{12.0}{rm}{\color[rgb]{0,0,0}$w_3$}%
}}}
\end{picture}

%% file: EPS-A/redge2ci1.pstex_t
\begin{picture}(0,0)%
\includegraphics{redge2ci1.pstex}%
\end{picture}%
\setlength{\unitlength}{3158sp}%
\begingroup\makeatletter\ifx\SetFigFont\undefined
\def\x#1#2#3#4#5#6#7\relax{\def\x{#1#2#3#4#5#6}}%
\expandafter\x\fmtname xxxxxx\relax \def\y{splain}%
\ifx\x\y   
\gdef\SetFigFont#1#2#3{%
  \ifnum #1<17\tiny\else \ifnum #1<20\small\else
  \ifnum #1<24\normalsize\else \ifnum #1<29\large\else
  \ifnum #1<34\Large\else \ifnum #1<41\LARGE\else
     \huge\fi\fi\fi\fi\fi\fi
  \csname #3\endcsname}%
\else
\gdef\SetFigFont#1#2#3{\begingroup
  \count@#1\relax \ifnum 25<\count@\count@25\fi
  \def\x{\endgroup\@setsize\SetFigFont{#2pt}}%
  \expandafter\x
    \csname \romannumeral\the\count@ pt\expandafter\endcsname
    \csname @\romannumeral\the\count@ pt\endcsname
  \csname #3\endcsname}%
\fi
\fi\endgroup
\begin{picture}(2911,1545)(173,-1107)
\put(526,-136){\makebox(0,0)[lb]{\smash{\SetFigFont{10}{12.0}{rm}{\color[rgb]{0,0,0}$v$}%
}}}
\put(908,-406){\makebox(0,0)[lb]{\smash{\SetFigFont{10}{12.0}{rm}{\color[rgb]{0,0,0}$1/2$}%
}}}
\put(1276,306){\makebox(0,0)[lb]{\smash{\SetFigFont{10}{12.0}{rm}{\color[rgb]{0,0,0}$1/2$}%
}}}
\put(1756,-398){\makebox(0,0)[lb]{\smash{\SetFigFont{10}{12.0}{rm}{\color[rgb]{0,0,0}$1/2$}%
}}}
\put(2446,-398){\makebox(0,0)[lb]{\smash{\SetFigFont{10}{12.0}{rm}{\color[rgb]{0,0,0}$1$}%
}}}
\put(1861,-1058){\makebox(0,0)[lb]{\smash{\SetFigFont{10}{12.0}{rm}{\color[rgb]{0,0,0}$c=5/2$}%
}}}
\put(3008,-83){\makebox(0,0)[lb]{\smash{\SetFigFont{10}{12.0}{rm}{\color[rgb]{0,0,0}$v'$}%
}}}
\put(1254,-83){\makebox(0,0)[lb]{\smash{\SetFigFont{10}{12.0}{rm}{\color[rgb]{0,0,0}$w_1$}%
}}}
\put(2221,-52){\makebox(0,0)[lb]{\smash{\SetFigFont{10}{12.0}{rm}{\color[rgb]{0,0,0}$w_2$}%
}}}
\end{picture}

%% file: EPS-A/redge2ci2.pstex_t
\begin{picture}(0,0)%
\includegraphics{redge2ci2.pstex}%
\end{picture}%
\setlength{\unitlength}{2763sp}%
\begingroup\makeatletter\ifx\SetFigFont\undefined
\def\x#1#2#3#4#5#6#7\relax{\def\x{#1#2#3#4#5#6}}%
\expandafter\x\fmtname xxxxxx\relax \def\y{splain}%
\ifx\x\y   
\gdef\SetFigFont#1#2#3{%
  \ifnum #1<17\tiny\else \ifnum #1<20\small\else
  \ifnum #1<24\normalsize\else \ifnum #1<29\large\else
  \ifnum #1<34\Large\else \ifnum #1<41\LARGE\else
     \huge\fi\fi\fi\fi\fi\fi
  \csname #3\endcsname}%
\else
\gdef\SetFigFont#1#2#3{\begingroup
  \count@#1\relax \ifnum 25<\count@\count@25\fi
  \def\x{\endgroup\@setsize\SetFigFont{#2pt}}%
  \expandafter\x
    \csname \romannumeral\the\count@ pt\expandafter\endcsname
    \csname @\romannumeral\the\count@ pt\endcsname
  \csname #3\endcsname}%
\fi
\fi\endgroup
\begin{picture}(2911,1545)(173,-1107)
\put(526,-136){\makebox(0,0)[lb]{\smash{\SetFigFont{8}{9.6}{rm}{\color[rgb]{0,0,0}$v$}%
}}}
\put(908,-406){\makebox(0,0)[lb]{\smash{\SetFigFont{8}{9.6}{rm}{\color[rgb]{0,0,0}$1/2$}%
}}}
\put(1276,306){\makebox(0,0)[lb]{\smash{\SetFigFont{8}{9.6}{rm}{\color[rgb]{0,0,0}$1/2$}%
}}}
\put(2266,-705){\makebox(0,0)[lb]{\smash{\SetFigFont{8}{9.6}{rm}{\color[rgb]{0,0,0}$1/2$}%
}}}
\put(1156,-1058){\makebox(0,0)[lb]{\smash{\SetFigFont{8}{9.6}{rm}{\color[rgb]{0,0,0}$c=3/2$}%
}}}
\put(1254,-83){\makebox(0,0)[lb]{\smash{\SetFigFont{8}{9.6}{rm}{\color[rgb]{0,0,0}$w_1$}%
}}}
\put(2221,-52){\makebox(0,0)[lb]{\smash{\SetFigFont{8}{9.6}{rm}{\color[rgb]{0,0,0}$w_2$}%
}}}
\put(3008,-83){\makebox(0,0)[lb]{\smash{\SetFigFont{8}{9.6}{rm}{\color[rgb]{0,0,0}$w_3$}%
}}}
\end{picture}

%% file: EPS-A/redge2cii.pstex_t
\begin{picture}(0,0)%
\includegraphics{redge2cii.pstex}%
\end{picture}%
\setlength{\unitlength}{3158sp}%
\begingroup\makeatletter\ifx\SetFigFont\undefined
\def\x#1#2#3#4#5#6#7\relax{\def\x{#1#2#3#4#5#6}}%
\expandafter\x\fmtname xxxxxx\relax \def\y{splain}%
\ifx\x\y   
\gdef\SetFigFont#1#2#3{%
  \ifnum #1<17\tiny\else \ifnum #1<20\small\else
  \ifnum #1<24\normalsize\else \ifnum #1<29\large\else
  \ifnum #1<34\Large\else \ifnum #1<41\LARGE\else
     \huge\fi\fi\fi\fi\fi\fi
  \csname #3\endcsname}%
\else
\gdef\SetFigFont#1#2#3{\begingroup
  \count@#1\relax \ifnum 25<\count@\count@25\fi
  \def\x{\endgroup\@setsize\SetFigFont{#2pt}}%
  \expandafter\x
    \csname \romannumeral\the\count@ pt\expandafter\endcsname
    \csname @\romannumeral\the\count@ pt\endcsname
  \csname #3\endcsname}%
\fi
\fi\endgroup
\begin{picture}(3308,1545)(526,-1107)
\put(526,-136){\makebox(0,0)[lb]{\smash{\SetFigFont{10}{12.0}{rm}{\color[rgb]{0,0,0}$v$}%
}}}
\put(908,-406){\makebox(0,0)[lb]{\smash{\SetFigFont{10}{12.0}{rm}{\color[rgb]{0,0,0}$1/2$}%
}}}
\put(1276,306){\makebox(0,0)[lb]{\smash{\SetFigFont{10}{12.0}{rm}{\color[rgb]{0,0,0}$1/2$}%
}}}
\put(1156,-1058){\makebox(0,0)[lb]{\smash{\SetFigFont{10}{12.0}{rm}{\color[rgb]{0,0,0}$c=3/2$}%
}}}
\put(2851,-930){\makebox(0,0)[lb]{\smash{\SetFigFont{10}{12.0}{rm}{\color[rgb]{0,0,0}$1/4$}%
}}}
\put(2348,-600){\makebox(0,0)[lb]{\smash{\SetFigFont{10}{12.0}{rm}{\color[rgb]{0,0,0}$1/4$}%
}}}
\put(1254,-83){\makebox(0,0)[lb]{\smash{\SetFigFont{10}{12.0}{rm}{\color[rgb]{0,0,0}$w_1$}%
}}}
\put(2221,-52){\makebox(0,0)[lb]{\smash{\SetFigFont{10}{12.0}{rm}{\color[rgb]{0,0,0}$w_2$}%
}}}
\put(3008,-83){\makebox(0,0)[lb]{\smash{\SetFigFont{10}{12.0}{rm}{\color[rgb]{0,0,0}$w_3$}%
}}}
\put(3766,-98){\makebox(0,0)[lb]{\smash{\SetFigFont{10}{12.0}{rm}{\color[rgb]{0,0,0}$w_4$}%
}}}
\end{picture}

%% file: EPS-A/redge2cii1.pstex_t
\begin{picture}(0,0)%
\includegraphics{redge2cii1.pstex}%
\end{picture}%
\setlength{\unitlength}{3158sp}%
\begingroup\makeatletter\ifx\SetFigFont\undefined
\def\x#1#2#3#4#5#6#7\relax{\def\x{#1#2#3#4#5#6}}%
\expandafter\x\fmtname xxxxxx\relax \def\y{splain}%
\ifx\x\y   
\gdef\SetFigFont#1#2#3{%
  \ifnum #1<17\tiny\else \ifnum #1<20\small\else
  \ifnum #1<24\normalsize\else \ifnum #1<29\large\else
  \ifnum #1<34\Large\else \ifnum #1<41\LARGE\else
     \huge\fi\fi\fi\fi\fi\fi
  \csname #3\endcsname}%
\else
\gdef\SetFigFont#1#2#3{\begingroup
  \count@#1\relax \ifnum 25<\count@\count@25\fi
  \def\x{\endgroup\@setsize\SetFigFont{#2pt}}%
  \expandafter\x
    \csname \romannumeral\the\count@ pt\expandafter\endcsname
    \csname @\romannumeral\the\count@ pt\endcsname
  \csname #3\endcsname}%
\fi
\fi\endgroup
\begin{picture}(3435,1596)(481,-1039)
\put(526,-136){\makebox(0,0)[lb]{\smash{\SetFigFont{10}{12.0}{rm}{\color[rgb]{0,0,0}$v$}%
}}}
\put(908,-406){\makebox(0,0)[lb]{\smash{\SetFigFont{10}{12.0}{rm}{\color[rgb]{0,0,0}$1/2$}%
}}}
\put(1276,306){\makebox(0,0)[lb]{\smash{\SetFigFont{10}{12.0}{rm}{\color[rgb]{0,0,0}$1/2$}%
}}}
\put(481,-990){\makebox(0,0)[lb]{\smash{\SetFigFont{10}{12.0}{rm}{\color[rgb]{0,0,0}$c=5/2$}%
}}}
\put(2266,-705){\makebox(0,0)[lb]{\smash{\SetFigFont{10}{12.0}{rm}{\color[rgb]{0,0,0}$1/4$}%
}}}
\put(2866,-975){\makebox(0,0)[lb]{\smash{\SetFigFont{10}{12.0}{rm}{\color[rgb]{0,0,0}$1/4$}%
}}}
\put(2619,382){\makebox(0,0)[lb]{\smash{\SetFigFont{10}{12.0}{rm}{\color[rgb]{0,0,0}$1/4$}%
}}}
\put(2438,-158){\makebox(0,0)[lb]{\smash{\SetFigFont{10}{12.0}{rm}{\color[rgb]{0,0,0}$1/4$}%
}}}
\put(3061,104){\makebox(0,0)[lb]{\smash{\SetFigFont{10}{12.0}{rm}{\color[rgb]{0,0,0}$1/2$}%
}}}
\put(3916,404){\makebox(0,0)[lb]{\smash{\SetFigFont{10}{12.0}{rm}{\color[rgb]{0,0,0}$x$}%
}}}
\put(3915, 82){\makebox(0,0)[lb]{\smash{\SetFigFont{10}{12.0}{rm}{\color[rgb]{0,0,0}$y$}%
}}}
\put(1254,-83){\makebox(0,0)[lb]{\smash{\SetFigFont{10}{12.0}{rm}{\color[rgb]{0,0,0}$w_1$}%
}}}
\put(3105,-271){\makebox(0,0)[lb]{\smash{\SetFigFont{10}{12.0}{rm}{\color[rgb]{0,0,0}$w_3$}%
}}}
\put(3893,-323){\makebox(0,0)[lb]{\smash{\SetFigFont{10}{12.0}{rm}{\color[rgb]{0,0,0}$w_4$}%
}}}
\put(1876,-361){\makebox(0,0)[lb]{\smash{\SetFigFont{10}{12.0}{rm}{\color[rgb]{0,0,0}$w_2$}%
}}}
\end{picture}

%% file: EPS-A/redge2cii2.pstex_t
\begin{picture}(0,0)%
\includegraphics{redge2cii2.pstex}%
\end{picture}%
\setlength{\unitlength}{3158sp}%
\begingroup\makeatletter\ifx\SetFigFont\undefined
\def\x#1#2#3#4#5#6#7\relax{\def\x{#1#2#3#4#5#6}}%
\expandafter\x\fmtname xxxxxx\relax \def\y{splain}%
\ifx\x\y   
\gdef\SetFigFont#1#2#3{%
  \ifnum #1<17\tiny\else \ifnum #1<20\small\else
  \ifnum #1<24\normalsize\else \ifnum #1<29\large\else
  \ifnum #1<34\Large\else \ifnum #1<41\LARGE\else
     \huge\fi\fi\fi\fi\fi\fi
  \csname #3\endcsname}%
\else
\gdef\SetFigFont#1#2#3{\begingroup
  \count@#1\relax \ifnum 25<\count@\count@25\fi
  \def\x{\endgroup\@setsize\SetFigFont{#2pt}}%
  \expandafter\x
    \csname \romannumeral\the\count@ pt\expandafter\endcsname
    \csname @\romannumeral\the\count@ pt\endcsname
  \csname #3\endcsname}%
\fi
\fi\endgroup
\begin{picture}(3427,1916)(481,-1039)
\put(526,-136){\makebox(0,0)[lb]{\smash{\SetFigFont{10}{12.0}{rm}{\color[rgb]{0,0,0}$v$}%
}}}
\put(908,-406){\makebox(0,0)[lb]{\smash{\SetFigFont{10}{12.0}{rm}{\color[rgb]{0,0,0}$1/2$}%
}}}
\put(1276,306){\makebox(0,0)[lb]{\smash{\SetFigFont{10}{12.0}{rm}{\color[rgb]{0,0,0}$1/2$}%
}}}
\put(2266,-705){\makebox(0,0)[lb]{\smash{\SetFigFont{10}{12.0}{rm}{\color[rgb]{0,0,0}$1/4$}%
}}}
\put(2866,-975){\makebox(0,0)[lb]{\smash{\SetFigFont{10}{12.0}{rm}{\color[rgb]{0,0,0}$1/4$}%
}}}
\put(2784,719){\makebox(0,0)[lb]{\smash{\SetFigFont{10}{12.0}{rm}{\color[rgb]{0,0,0}$1/4$}%
}}}
\put(2912,397){\makebox(0,0)[lb]{\smash{\SetFigFont{10}{12.0}{rm}{\color[rgb]{0,0,0}$1/4$}%
}}}
\put(3894,719){\makebox(0,0)[lb]{\smash{\SetFigFont{10}{12.0}{rm}{\color[rgb]{0,0,0}$x$}%
}}}
\put(3908,412){\makebox(0,0)[lb]{\smash{\SetFigFont{10}{12.0}{rm}{\color[rgb]{0,0,0}$y$}%
}}}
\put(3901, 82){\makebox(0,0)[lb]{\smash{\SetFigFont{10}{12.0}{rm}{\color[rgb]{0,0,0}$z$}%
}}}
\put(3047, 82){\makebox(0,0)[lb]{\smash{\SetFigFont{10}{12.0}{rm}{\color[rgb]{0,0,0}$1/4$}%
}}}
\put(481,-990){\makebox(0,0)[lb]{\smash{\SetFigFont{10}{12.0}{rm}{\color[rgb]{0,0,0}$c=9/4$}%
}}}
\put(1254,-83){\makebox(0,0)[lb]{\smash{\SetFigFont{10}{12.0}{rm}{\color[rgb]{0,0,0}$w_1$}%
}}}
\put(1906,-262){\makebox(0,0)[lb]{\smash{\SetFigFont{10}{12.0}{rm}{\color[rgb]{0,0,0}$w_2$}%
}}}
\put(3105,-271){\makebox(0,0)[lb]{\smash{\SetFigFont{10}{12.0}{rm}{\color[rgb]{0,0,0}$w_3$}%
}}}
\put(3893,-323){\makebox(0,0)[lb]{\smash{\SetFigFont{10}{12.0}{rm}{\color[rgb]{0,0,0}$w_4$}%
}}}
\end{picture}

%% file: EPS-A/redge2cii3.pstex_t
\begin{picture}(0,0)%
\includegraphics{redge2cii3.pstex}%
\end{picture}%
\setlength{\unitlength}{3158sp}%
\begingroup\makeatletter\ifx\SetFigFont\undefined
\def\x#1#2#3#4#5#6#7\relax{\def\x{#1#2#3#4#5#6}}%
\expandafter\x\fmtname xxxxxx\relax \def\y{splain}%
\ifx\x\y   
\gdef\SetFigFont#1#2#3{%
  \ifnum #1<17\tiny\else \ifnum #1<20\small\else
  \ifnum #1<24\normalsize\else \ifnum #1<29\large\else
  \ifnum #1<34\Large\else \ifnum #1<41\LARGE\else
     \huge\fi\fi\fi\fi\fi\fi
  \csname #3\endcsname}%
\else
\gdef\SetFigFont#1#2#3{\begingroup
  \count@#1\relax \ifnum 25<\count@\count@25\fi
  \def\x{\endgroup\@setsize\SetFigFont{#2pt}}%
  \expandafter\x
    \csname \romannumeral\the\count@ pt\expandafter\endcsname
    \csname @\romannumeral\the\count@ pt\endcsname
  \csname #3\endcsname}%
\fi
\fi\endgroup
\begin{picture}(3420,2235)(331,-1189)
\put(376,-286){\makebox(0,0)[lb]{\smash{\SetFigFont{10}{12.0}{rm}{\color[rgb]{0,0,0}$v$}%
}}}
\put(758,-556){\makebox(0,0)[lb]{\smash{\SetFigFont{10}{12.0}{rm}{\color[rgb]{0,0,0}$1/2$}%
}}}
\put(1126,156){\makebox(0,0)[lb]{\smash{\SetFigFont{10}{12.0}{rm}{\color[rgb]{0,0,0}$1/2$}%
}}}
\put(2116,-855){\makebox(0,0)[lb]{\smash{\SetFigFont{10}{12.0}{rm}{\color[rgb]{0,0,0}$1/4$}%
}}}
\put(2716,-1125){\makebox(0,0)[lb]{\smash{\SetFigFont{10}{12.0}{rm}{\color[rgb]{0,0,0}$1/4$}%
}}}
\put(2326,689){\makebox(0,0)[lb]{\smash{\SetFigFont{10}{12.0}{rm}{\color[rgb]{0,0,0}$1/4$}%
}}}
\put(2776,614){\makebox(0,0)[lb]{\smash{\SetFigFont{10}{12.0}{rm}{\color[rgb]{0,0,0}$1/4$}%
}}}
\put(2776, 89){\makebox(0,0)[lb]{\smash{\SetFigFont{10}{12.0}{rm}{\color[rgb]{0,0,0}$1/8$}%
}}}
\put(3151,-286){\makebox(0,0)[lb]{\smash{\SetFigFont{10}{12.0}{rm}{\color[rgb]{0,0,0}$1/8$}%
}}}
\put(3751,914){\makebox(0,0)[lb]{\smash{\SetFigFont{10}{12.0}{rm}{\color[rgb]{0,0,0}$x$}%
}}}
\put(3751,614){\makebox(0,0)[lb]{\smash{\SetFigFont{10}{12.0}{rm}{\color[rgb]{0,0,0}$y$}%
}}}
\put(3751,314){\makebox(0,0)[lb]{\smash{\SetFigFont{10}{12.0}{rm}{\color[rgb]{0,0,0}$z$}%
}}}
\put(3751, 14){\makebox(0,0)[lb]{\smash{\SetFigFont{10}{12.0}{rm}{\color[rgb]{0,0,0}$t$}%
}}}
\put(331,-1140){\makebox(0,0)[lb]{\smash{\SetFigFont{10}{12.0}{rm}{\color[rgb]{0,0,0}$c=9/4$}%
}}}
\put(1104,-233){\makebox(0,0)[lb]{\smash{\SetFigFont{10}{12.0}{rm}{\color[rgb]{0,0,0}$w_1$}%
}}}
\put(1756,-412){\makebox(0,0)[lb]{\smash{\SetFigFont{10}{12.0}{rm}{\color[rgb]{0,0,0}$w_2$}%
}}}
\put(2776,-586){\makebox(0,0)[lb]{\smash{\SetFigFont{10}{12.0}{rm}{\color[rgb]{0,0,0}$w_3$}%
}}}
\put(3743,-473){\makebox(0,0)[lb]{\smash{\SetFigFont{10}{12.0}{rm}{\color[rgb]{0,0,0}$w_4$}%
}}}
\end{picture}

%% file: EPS-A/redge2cii4.pstex_t
\begin{picture}(0,0)%
\includegraphics{redge2cii4.pstex}%
\end{picture}%
\setlength{\unitlength}{3158sp}%
\begingroup\makeatletter\ifx\SetFigFont\undefined
\def\x#1#2#3#4#5#6#7\relax{\def\x{#1#2#3#4#5#6}}%
\expandafter\x\fmtname xxxxxx\relax \def\y{splain}%
\ifx\x\y   
\gdef\SetFigFont#1#2#3{%
  \ifnum #1<17\tiny\else \ifnum #1<20\small\else
  \ifnum #1<24\normalsize\else \ifnum #1<29\large\else
  \ifnum #1<34\Large\else \ifnum #1<41\LARGE\else
     \huge\fi\fi\fi\fi\fi\fi
  \csname #3\endcsname}%
\else
\gdef\SetFigFont#1#2#3{\begingroup
  \count@#1\relax \ifnum 25<\count@\count@25\fi
  \def\x{\endgroup\@setsize\SetFigFont{#2pt}}%
  \expandafter\x
    \csname \romannumeral\the\count@ pt\expandafter\endcsname
    \csname @\romannumeral\the\count@ pt\endcsname
  \csname #3\endcsname}%
\fi
\fi\endgroup
\begin{picture}(4170,2235)(331,-1189)
\put(376,-286){\makebox(0,0)[lb]{\smash{\SetFigFont{10}{12.0}{rm}{\color[rgb]{0,0,0}$v$}%
}}}
\put(758,-556){\makebox(0,0)[lb]{\smash{\SetFigFont{10}{12.0}{rm}{\color[rgb]{0,0,0}$1/2$}%
}}}
\put(1126,156){\makebox(0,0)[lb]{\smash{\SetFigFont{10}{12.0}{rm}{\color[rgb]{0,0,0}$1/2$}%
}}}
\put(2716,-1125){\makebox(0,0)[lb]{\smash{\SetFigFont{10}{12.0}{rm}{\color[rgb]{0,0,0}$1/4$}%
}}}
\put(2326,689){\makebox(0,0)[lb]{\smash{\SetFigFont{10}{12.0}{rm}{\color[rgb]{0,0,0}$1/4$}%
}}}
\put(2776, 89){\makebox(0,0)[lb]{\smash{\SetFigFont{10}{12.0}{rm}{\color[rgb]{0,0,0}$1/8$}%
}}}
\put(3151,-286){\makebox(0,0)[lb]{\smash{\SetFigFont{10}{12.0}{rm}{\color[rgb]{0,0,0}$1/8$}%
}}}
\put(3751,914){\makebox(0,0)[lb]{\smash{\SetFigFont{10}{12.0}{rm}{\color[rgb]{0,0,0}$x$}%
}}}
\put(3751,314){\makebox(0,0)[lb]{\smash{\SetFigFont{10}{12.0}{rm}{\color[rgb]{0,0,0}$z$}%
}}}
\put(3751, 14){\makebox(0,0)[lb]{\smash{\SetFigFont{10}{12.0}{rm}{\color[rgb]{0,0,0}$t$}%
}}}
\put(331,-1140){\makebox(0,0)[lb]{\smash{\SetFigFont{10}{12.0}{rm}{\color[rgb]{0,0,0}$c=11/4$}%
}}}
\put(4501,-361){\makebox(0,0)[lb]{\smash{\SetFigFont{10}{12.0}{rm}{\color[rgb]{0,0,0}$w$}%
}}}
\put(2701,-886){\makebox(0,0)[lb]{\smash{\SetFigFont{10}{12.0}{rm}{\color[rgb]{0,0,0}$1/4$}%
}}}
\put(1876,-811){\makebox(0,0)[lb]{\smash{\SetFigFont{10}{12.0}{rm}{\color[rgb]{0,0,0}$1/4$}%
}}}
\put(3826,-286){\makebox(0,0)[lb]{\smash{\SetFigFont{10}{12.0}{rm}{\color[rgb]{0,0,0}$1/2$}%
}}}
\put(1104,-233){\makebox(0,0)[lb]{\smash{\SetFigFont{10}{12.0}{rm}{\color[rgb]{0,0,0}$w_1$}%
}}}
\put(1756,-412){\makebox(0,0)[lb]{\smash{\SetFigFont{10}{12.0}{rm}{\color[rgb]{0,0,0}$w_2$}%
}}}
\put(2776,-586){\makebox(0,0)[lb]{\smash{\SetFigFont{10}{12.0}{rm}{\color[rgb]{0,0,0}$w_3$}%
}}}
\put(3676,-586){\makebox(0,0)[lb]{\smash{\SetFigFont{10}{12.0}{rm}{\color[rgb]{0,0,0}$w_4$}%
}}}
\end{picture}

%% file: EPS-A/redge2ciii.pstex_t
\begin{picture}(0,0)%
\includegraphics{redge2ciii.pstex}%
\end{picture}%
\setlength{\unitlength}{3158sp}%
\begingroup\makeatletter\ifx\SetFigFont\undefined
\def\x#1#2#3#4#5#6#7\relax{\def\x{#1#2#3#4#5#6}}%
\expandafter\x\fmtname xxxxxx\relax \def\y{splain}%
\ifx\x\y   
\gdef\SetFigFont#1#2#3{%
  \ifnum #1<17\tiny\else \ifnum #1<20\small\else
  \ifnum #1<24\normalsize\else \ifnum #1<29\large\else
  \ifnum #1<34\Large\else \ifnum #1<41\LARGE\else
     \huge\fi\fi\fi\fi\fi\fi
  \csname #3\endcsname}%
\else
\gdef\SetFigFont#1#2#3{\begingroup
  \count@#1\relax \ifnum 25<\count@\count@25\fi
  \def\x{\endgroup\@setsize\SetFigFont{#2pt}}%
  \expandafter\x
    \csname \romannumeral\the\count@ pt\expandafter\endcsname
    \csname @\romannumeral\the\count@ pt\endcsname
  \csname #3\endcsname}%
\fi
\fi\endgroup
\begin{picture}(2648,1477)(481,-1039)
\put(526,-136){\makebox(0,0)[lb]{\smash{\SetFigFont{10}{12.0}{rm}{\color[rgb]{0,0,0}$v$}%
}}}
\put(908,-406){\makebox(0,0)[lb]{\smash{\SetFigFont{10}{12.0}{rm}{\color[rgb]{0,0,0}$1/2$}%
}}}
\put(1276,306){\makebox(0,0)[lb]{\smash{\SetFigFont{10}{12.0}{rm}{\color[rgb]{0,0,0}$1/2$}%
}}}
\put(2266,-705){\makebox(0,0)[lb]{\smash{\SetFigFont{10}{12.0}{rm}{\color[rgb]{0,0,0}$1/4$}%
}}}
\put(481,-990){\makebox(0,0)[lb]{\smash{\SetFigFont{10}{12.0}{rm}{\color[rgb]{0,0,0}$c=9/4$}%
}}}
\put(3129, 89){\makebox(0,0)[lb]{\smash{\SetFigFont{10}{12.0}{rm}{\color[rgb]{0,0,0}$x$}%
}}}
\put(1718,-375){\makebox(0,0)[lb]{\smash{\SetFigFont{10}{12.0}{rm}{\color[rgb]{0,0,0}$1/4$}%
}}}
\put(2356,120){\makebox(0,0)[lb]{\smash{\SetFigFont{10}{12.0}{rm}{\color[rgb]{0,0,0}$3/4$}%
}}}
\put(1254,-83){\makebox(0,0)[lb]{\smash{\SetFigFont{10}{12.0}{rm}{\color[rgb]{0,0,0}$w_1$}%
}}}
\put(3113,-338){\makebox(0,0)[lb]{\smash{\SetFigFont{10}{12.0}{rm}{\color[rgb]{0,0,0}$w_4$}%
}}}
\put(2364,-232){\makebox(0,0)[lb]{\smash{\SetFigFont{10}{12.0}{rm}{\color[rgb]{0,0,0}$w_2$}%
}}}
\end{picture}

%% file: EPS-A/redge2ciiix.pstex_t
\begin{picture}(0,0)%
\includegraphics{redge2ciiix.pstex}%
\end{picture}%
\setlength{\unitlength}{3158sp}%
\begingroup\makeatletter\ifx\SetFigFont\undefined
\def\x#1#2#3#4#5#6#7\relax{\def\x{#1#2#3#4#5#6}}%
\expandafter\x\fmtname xxxxxx\relax \def\y{splain}%
\ifx\x\y   
\gdef\SetFigFont#1#2#3{%
  \ifnum #1<17\tiny\else \ifnum #1<20\small\else
  \ifnum #1<24\normalsize\else \ifnum #1<29\large\else
  \ifnum #1<34\Large\else \ifnum #1<41\LARGE\else
     \huge\fi\fi\fi\fi\fi\fi
  \csname #3\endcsname}%
\else
\gdef\SetFigFont#1#2#3{\begingroup
  \count@#1\relax \ifnum 25<\count@\count@25\fi
  \def\x{\endgroup\@setsize\SetFigFont{#2pt}}%
  \expandafter\x
    \csname \romannumeral\the\count@ pt\expandafter\endcsname
    \csname @\romannumeral\the\count@ pt\endcsname
  \csname #3\endcsname}%
\fi
\fi\endgroup
\begin{picture}(3789,761)(-26,-542)
\put(503,-106){\makebox(0,0)[lb]{\smash{\SetFigFont{10}{12.0}{rm}{\color[rgb]{0,0,0}$v$}%
}}}
\put(1254,-83){\makebox(0,0)[lb]{\smash{\SetFigFont{10}{12.0}{rm}{\color[rgb]{0,0,0}$w_1$}%
}}}
\put(2207,-68){\makebox(0,0)[lb]{\smash{\SetFigFont{10}{12.0}{rm}{\color[rgb]{0,0,0}$w_2$}%
}}}
\put(3060,-413){\makebox(0,0)[lb]{\smash{\SetFigFont{10}{12.0}{rm}{\color[rgb]{0,0,0}$w_4$}%
}}}
\end{picture}

%% file: EPS-A/redge2ciii1a.pstex_t
\begin{picture}(0,0)%
\includegraphics{redge2ciii1a.pstex}%
\end{picture}%
\setlength{\unitlength}{3158sp}%
\begingroup\makeatletter\ifx\SetFigFont\undefined
\def\x#1#2#3#4#5#6#7\relax{\def\x{#1#2#3#4#5#6}}%
\expandafter\x\fmtname xxxxxx\relax \def\y{splain}%
\ifx\x\y   
\gdef\SetFigFont#1#2#3{%
  \ifnum #1<17\tiny\else \ifnum #1<20\small\else
  \ifnum #1<24\normalsize\else \ifnum #1<29\large\else
  \ifnum #1<34\Large\else \ifnum #1<41\LARGE\else
     \huge\fi\fi\fi\fi\fi\fi
  \csname #3\endcsname}%
\else
\gdef\SetFigFont#1#2#3{\begingroup
  \count@#1\relax \ifnum 25<\count@\count@25\fi
  \def\x{\endgroup\@setsize\SetFigFont{#2pt}}%
  \expandafter\x
    \csname \romannumeral\the\count@ pt\expandafter\endcsname
    \csname @\romannumeral\the\count@ pt\endcsname
  \csname #3\endcsname}%
\fi
\fi\endgroup
\begin{picture}(4118,1847)(503,-1181)
\put(526,-136){\makebox(0,0)[lb]{\smash{\SetFigFont{10}{12.0}{rm}{\color[rgb]{0,0,0}$v$}%
}}}
\put(908,-406){\makebox(0,0)[lb]{\smash{\SetFigFont{10}{12.0}{rm}{\color[rgb]{0,0,0}$1/2$}%
}}}
\put(3061,-98){\makebox(0,0)[lb]{\smash{\SetFigFont{10}{12.0}{rm}{\color[rgb]{0,0,0}$x$}%
}}}
\put(1778,-375){\makebox(0,0)[lb]{\smash{\SetFigFont{10}{12.0}{rm}{\color[rgb]{0,0,0}$1/4$}%
}}}
\put(2723,-781){\makebox(0,0)[lb]{\smash{\SetFigFont{10}{12.0}{rm}{\color[rgb]{0,0,0}$1/4$}%
}}}
\put(2408,-382){\makebox(0,0)[lb]{\smash{\SetFigFont{10}{12.0}{rm}{\color[rgb]{0,0,0}$3/4$}%
}}}
\put(3120,-367){\makebox(0,0)[lb]{\smash{\SetFigFont{10}{12.0}{rm}{\color[rgb]{0,0,0}$3/4$}%
}}}
\put(3998,-383){\makebox(0,0)[lb]{\smash{\SetFigFont{10}{12.0}{rm}{\color[rgb]{0,0,0}$1$}%
}}}
\put(1344,261){\makebox(0,0)[lb]{\smash{\SetFigFont{10}{12.0}{rm}{\color[rgb]{0,0,0}$1/2$}%
}}}
\put(4621,-151){\makebox(0,0)[lb]{\smash{\SetFigFont{10}{12.0}{rm}{\color[rgb]{0,0,0}$z$}%
}}}
\put(1171,-1132){\makebox(0,0)[lb]{\smash{\SetFigFont{10}{12.0}{rm}{\color[rgb]{0,0,0}$c=4$}%
}}}
\put(1254,-83){\makebox(0,0)[lb]{\smash{\SetFigFont{10}{12.0}{rm}{\color[rgb]{0,0,0}$w_1$}%
}}}
\put(2274,-105){\makebox(0,0)[lb]{\smash{\SetFigFont{10}{12.0}{rm}{\color[rgb]{0,0,0}$w_2$}%
}}}
\put(3736,-98){\makebox(0,0)[lb]{\smash{\SetFigFont{10}{12.0}{rm}{\color[rgb]{0,0,0}$w_4$}%
}}}
\end{picture}

%% file: EPS-A/redge2ciii1b.pstex_t
\begin{picture}(0,0)%
\includegraphics{redge2ciii1b.pstex}%
\end{picture}%
\setlength{\unitlength}{3158sp}%
\begingroup\makeatletter\ifx\SetFigFont\undefined
\def\x#1#2#3#4#5#6#7\relax{\def\x{#1#2#3#4#5#6}}%
\expandafter\x\fmtname xxxxxx\relax \def\y{splain}%
\ifx\x\y   
\gdef\SetFigFont#1#2#3{%
  \ifnum #1<17\tiny\else \ifnum #1<20\small\else
  \ifnum #1<24\normalsize\else \ifnum #1<29\large\else
  \ifnum #1<34\Large\else \ifnum #1<41\LARGE\else
     \huge\fi\fi\fi\fi\fi\fi
  \csname #3\endcsname}%
\else
\gdef\SetFigFont#1#2#3{\begingroup
  \count@#1\relax \ifnum 25<\count@\count@25\fi
  \def\x{\endgroup\@setsize\SetFigFont{#2pt}}%
  \expandafter\x
    \csname \romannumeral\the\count@ pt\expandafter\endcsname
    \csname @\romannumeral\the\count@ pt\endcsname
  \csname #3\endcsname}%
\fi
\fi\endgroup
\begin{picture}(4088,1847)(503,-1181)
\put(526,-136){\makebox(0,0)[lb]{\smash{\SetFigFont{10}{12.0}{rm}{\color[rgb]{0,0,0}$v$}%
}}}
\put(908,-406){\makebox(0,0)[lb]{\smash{\SetFigFont{10}{12.0}{rm}{\color[rgb]{0,0,0}$1/2$}%
}}}
\put(3061,-98){\makebox(0,0)[lb]{\smash{\SetFigFont{10}{12.0}{rm}{\color[rgb]{0,0,0}$x$}%
}}}
\put(1778,-375){\makebox(0,0)[lb]{\smash{\SetFigFont{10}{12.0}{rm}{\color[rgb]{0,0,0}$1/4$}%
}}}
\put(2408,-382){\makebox(0,0)[lb]{\smash{\SetFigFont{10}{12.0}{rm}{\color[rgb]{0,0,0}$3/4$}%
}}}
\put(3120,-367){\makebox(0,0)[lb]{\smash{\SetFigFont{10}{12.0}{rm}{\color[rgb]{0,0,0}$3/4$}%
}}}
\put(1344,261){\makebox(0,0)[lb]{\smash{\SetFigFont{10}{12.0}{rm}{\color[rgb]{0,0,0}$1/2$}%
}}}
\put(3706,-98){\makebox(0,0)[lb]{\smash{\SetFigFont{10}{12.0}{rm}{\color[rgb]{0,0,0}$y$}%
}}}
\put(2814,-638){\makebox(0,0)[lb]{\smash{\SetFigFont{10}{12.0}{rm}{\color[rgb]{0,0,0}$1/4$}%
}}}
\put(1171,-1132){\makebox(0,0)[lb]{\smash{\SetFigFont{10}{12.0}{rm}{\color[rgb]{0,0,0}$c=3$}%
}}}
\put(1254,-83){\makebox(0,0)[lb]{\smash{\SetFigFont{10}{12.0}{rm}{\color[rgb]{0,0,0}$w_1$}%
}}}
\put(2274,-105){\makebox(0,0)[lb]{\smash{\SetFigFont{10}{12.0}{rm}{\color[rgb]{0,0,0}$w_2$}%
}}}
\put(4538,-90){\makebox(0,0)[lb]{\smash{\SetFigFont{10}{12.0}{rm}{\color[rgb]{0,0,0}$w_4$}%
}}}
\end{picture}

%% file: EPS-A/redge2ciii2.pstex_t
\begin{picture}(0,0)%
\includegraphics{redge2ciii2.pstex}%
\end{picture}%
\setlength{\unitlength}{3158sp}%
\begingroup\makeatletter\ifx\SetFigFont\undefined
\def\x#1#2#3#4#5#6#7\relax{\def\x{#1#2#3#4#5#6}}%
\expandafter\x\fmtname xxxxxx\relax \def\y{splain}%
\ifx\x\y   
\gdef\SetFigFont#1#2#3{%
  \ifnum #1<17\tiny\else \ifnum #1<20\small\else
  \ifnum #1<24\normalsize\else \ifnum #1<29\large\else
  \ifnum #1<34\Large\else \ifnum #1<41\LARGE\else
     \huge\fi\fi\fi\fi\fi\fi
  \csname #3\endcsname}%
\else
\gdef\SetFigFont#1#2#3{\begingroup
  \count@#1\relax \ifnum 25<\count@\count@25\fi
  \def\x{\endgroup\@setsize\SetFigFont{#2pt}}%
  \expandafter\x
    \csname \romannumeral\the\count@ pt\expandafter\endcsname
    \csname @\romannumeral\the\count@ pt\endcsname
  \csname #3\endcsname}%
\fi
\fi\endgroup
\begin{picture}(4822,1574)(526,-1181)
\put(526,-136){\makebox(0,0)[lb]{\smash{\SetFigFont{10}{12.0}{rm}{\color[rgb]{0,0,0}$v$}%
}}}
\put(908,-406){\makebox(0,0)[lb]{\smash{\SetFigFont{10}{12.0}{rm}{\color[rgb]{0,0,0}$1/2$}%
}}}
\put(3061,-98){\makebox(0,0)[lb]{\smash{\SetFigFont{10}{12.0}{rm}{\color[rgb]{0,0,0}$x$}%
}}}
\put(1778,-375){\makebox(0,0)[lb]{\smash{\SetFigFont{10}{12.0}{rm}{\color[rgb]{0,0,0}$1/4$}%
}}}
\put(2408,-382){\makebox(0,0)[lb]{\smash{\SetFigFont{10}{12.0}{rm}{\color[rgb]{0,0,0}$3/4$}%
}}}
\put(1344,261){\makebox(0,0)[lb]{\smash{\SetFigFont{10}{12.0}{rm}{\color[rgb]{0,0,0}$1/2$}%
}}}
\put(3174,-608){\makebox(0,0)[lb]{\smash{\SetFigFont{10}{12.0}{rm}{\color[rgb]{0,0,0}$1/4$}%
}}}
\put(1171,-1132){\makebox(0,0)[lb]{\smash{\SetFigFont{10}{12.0}{rm}{\color[rgb]{0,0,0}$c=9/4$}%
}}}
\put(1254,-83){\makebox(0,0)[lb]{\smash{\SetFigFont{10}{12.0}{rm}{\color[rgb]{0,0,0}$w_1$}%
}}}
\put(2274,-105){\makebox(0,0)[lb]{\smash{\SetFigFont{10}{12.0}{rm}{\color[rgb]{0,0,0}$w_2$}%
}}}
\put(5348,-143){\makebox(0,0)[lb]{\smash{\SetFigFont{10}{12.0}{rm}{\color[rgb]{0,0,0}$w_4$}%
}}}
\end{picture}

%% file: EPS-A/redge2ciii2a.pstex_t
\begin{picture}(0,0)%
\includegraphics{redge2ciii2a.pstex}%
\end{picture}%
\setlength{\unitlength}{3158sp}%
\begingroup\makeatletter\ifx\SetFigFont\undefined
\def\x#1#2#3#4#5#6#7\relax{\def\x{#1#2#3#4#5#6}}%
\expandafter\x\fmtname xxxxxx\relax \def\y{splain}%
\ifx\x\y   
\gdef\SetFigFont#1#2#3{%
  \ifnum #1<17\tiny\else \ifnum #1<20\small\else
  \ifnum #1<24\normalsize\else \ifnum #1<29\large\else
  \ifnum #1<34\Large\else \ifnum #1<41\LARGE\else
     \huge\fi\fi\fi\fi\fi\fi
  \csname #3\endcsname}%
\else
\gdef\SetFigFont#1#2#3{\begingroup
  \count@#1\relax \ifnum 25<\count@\count@25\fi
  \def\x{\endgroup\@setsize\SetFigFont{#2pt}}%
  \expandafter\x
    \csname \romannumeral\the\count@ pt\expandafter\endcsname
    \csname @\romannumeral\the\count@ pt\endcsname
  \csname #3\endcsname}%
\fi
\fi\endgroup
\begin{picture}(4125,1574)(526,-1181)
\put(526,-136){\makebox(0,0)[lb]{\smash{\SetFigFont{10}{12.0}{rm}{\color[rgb]{0,0,0}$v$}%
}}}
\put(908,-406){\makebox(0,0)[lb]{\smash{\SetFigFont{10}{12.0}{rm}{\color[rgb]{0,0,0}$1/2$}%
}}}
\put(1778,-375){\makebox(0,0)[lb]{\smash{\SetFigFont{10}{12.0}{rm}{\color[rgb]{0,0,0}$1/4$}%
}}}
\put(2408,-382){\makebox(0,0)[lb]{\smash{\SetFigFont{10}{12.0}{rm}{\color[rgb]{0,0,0}$3/4$}%
}}}
\put(1344,261){\makebox(0,0)[lb]{\smash{\SetFigFont{10}{12.0}{rm}{\color[rgb]{0,0,0}$1/2$}%
}}}
\put(2926,-75){\makebox(0,0)[lb]{\smash{\SetFigFont{10}{12.0}{rm}{\color[rgb]{0,0,0}$x$}%
}}}
\put(4651,-286){\makebox(0,0)[lb]{\smash{\SetFigFont{10}{12.0}{rm}{\color[rgb]{0,0,0}$w$}%
}}}
\put(4651,104){\makebox(0,0)[lb]{\smash{\SetFigFont{10}{12.0}{rm}{\color[rgb]{0,0,0}$z$}%
}}}
\put(1171,-1132){\makebox(0,0)[lb]{\smash{\SetFigFont{10}{12.0}{rm}{\color[rgb]{0,0,0}$c=29/8$}%
}}}
\put(2709,-796){\makebox(0,0)[lb]{\smash{\SetFigFont{10}{12.0}{rm}{\color[rgb]{0,0,0}$1/4$}%
}}}
\put(3136,-369){\makebox(0,0)[lb]{\smash{\SetFigFont{10}{12.0}{rm}{\color[rgb]{0,0,0}$3/8$}%
}}}
\put(3414, 96){\makebox(0,0)[lb]{\smash{\SetFigFont{10}{12.0}{rm}{\color[rgb]{0,0,0}$3/8$}%
}}}
\put(3923,-188){\makebox(0,0)[lb]{\smash{\SetFigFont{10}{12.0}{rm}{\color[rgb]{0,0,0}$5/8$}%
}}}
\put(1254,-83){\makebox(0,0)[lb]{\smash{\SetFigFont{10}{12.0}{rm}{\color[rgb]{0,0,0}$w_1$}%
}}}
\put(2274,-105){\makebox(0,0)[lb]{\smash{\SetFigFont{10}{12.0}{rm}{\color[rgb]{0,0,0}$w_2$}%
}}}
\put(3810,-428){\makebox(0,0)[lb]{\smash{\SetFigFont{10}{12.0}{rm}{\color[rgb]{0,0,0}$w_4$}%
}}}
\end{picture}

%% file: EPS-A/redge2ciii2b.pstex_t
\begin{picture}(0,0)%
\includegraphics{redge2ciii2b.pstex}%
\end{picture}%
\setlength{\unitlength}{3158sp}%
\begingroup\makeatletter\ifx\SetFigFont\undefined
\def\x#1#2#3#4#5#6#7\relax{\def\x{#1#2#3#4#5#6}}%
\expandafter\x\fmtname xxxxxx\relax \def\y{splain}%
\ifx\x\y   
\gdef\SetFigFont#1#2#3{%
  \ifnum #1<17\tiny\else \ifnum #1<20\small\else
  \ifnum #1<24\normalsize\else \ifnum #1<29\large\else
  \ifnum #1<34\Large\else \ifnum #1<41\LARGE\else
     \huge\fi\fi\fi\fi\fi\fi
  \csname #3\endcsname}%
\else
\gdef\SetFigFont#1#2#3{\begingroup
  \count@#1\relax \ifnum 25<\count@\count@25\fi
  \def\x{\endgroup\@setsize\SetFigFont{#2pt}}%
  \expandafter\x
    \csname \romannumeral\the\count@ pt\expandafter\endcsname
    \csname @\romannumeral\the\count@ pt\endcsname
  \csname #3\endcsname}%
\fi
\fi\endgroup
\begin{picture}(3712,1803)(526,-1181)
\put(526,-136){\makebox(0,0)[lb]{\smash{\SetFigFont{10}{12.0}{rm}{\color[rgb]{0,0,0}$v$}%
}}}
\put(908,-406){\makebox(0,0)[lb]{\smash{\SetFigFont{10}{12.0}{rm}{\color[rgb]{0,0,0}$1/2$}%
}}}
\put(1778,-375){\makebox(0,0)[lb]{\smash{\SetFigFont{10}{12.0}{rm}{\color[rgb]{0,0,0}$1/4$}%
}}}
\put(2408,-382){\makebox(0,0)[lb]{\smash{\SetFigFont{10}{12.0}{rm}{\color[rgb]{0,0,0}$3/4$}%
}}}
\put(1344,261){\makebox(0,0)[lb]{\smash{\SetFigFont{10}{12.0}{rm}{\color[rgb]{0,0,0}$1/2$}%
}}}
\put(2709,-796){\makebox(0,0)[lb]{\smash{\SetFigFont{10}{12.0}{rm}{\color[rgb]{0,0,0}$1/4$}%
}}}
\put(4238, 89){\makebox(0,0)[lb]{\smash{\SetFigFont{10}{12.0}{rm}{\color[rgb]{0,0,0}$z$}%
}}}
\put(4230,471){\makebox(0,0)[lb]{\smash{\SetFigFont{10}{12.0}{rm}{\color[rgb]{0,0,0}$w$}%
}}}
\put(2813,-90){\makebox(0,0)[lb]{\smash{\SetFigFont{10}{12.0}{rm}{\color[rgb]{0,0,0}$x$}%
}}}
\put(3294,396){\makebox(0,0)[lb]{\smash{\SetFigFont{10}{12.0}{rm}{\color[rgb]{0,0,0}$3/8$}%
}}}
\put(3466,134){\makebox(0,0)[lb]{\smash{\SetFigFont{10}{12.0}{rm}{\color[rgb]{0,0,0}$3/8$}%
}}}
\put(1171,-1132){\makebox(0,0)[lb]{\smash{\SetFigFont{10}{12.0}{rm}{\color[rgb]{0,0,0}$c=3$}%
}}}
\put(1254,-83){\makebox(0,0)[lb]{\smash{\SetFigFont{10}{12.0}{rm}{\color[rgb]{0,0,0}$w_1$}%
}}}
\put(2274,-105){\makebox(0,0)[lb]{\smash{\SetFigFont{10}{12.0}{rm}{\color[rgb]{0,0,0}$w_2$}%
}}}
\put(3810,-428){\makebox(0,0)[lb]{\smash{\SetFigFont{10}{12.0}{rm}{\color[rgb]{0,0,0}$w_4$}%
}}}
\end{picture}

%% file: EPS-A/redge2civ1a.pstex_t
\begin{picture}(0,0)%
\includegraphics{redge2civ1a.pstex}%
\end{picture}%
\setlength{\unitlength}{3158sp}%
\begingroup\makeatletter\ifx\SetFigFont\undefined
\def\x#1#2#3#4#5#6#7\relax{\def\x{#1#2#3#4#5#6}}%
\expandafter\x\fmtname xxxxxx\relax \def\y{splain}%
\ifx\x\y   
\gdef\SetFigFont#1#2#3{%
  \ifnum #1<17\tiny\else \ifnum #1<20\small\else
  \ifnum #1<24\normalsize\else \ifnum #1<29\large\else
  \ifnum #1<34\Large\else \ifnum #1<41\LARGE\else
     \huge\fi\fi\fi\fi\fi\fi
  \csname #3\endcsname}%
\else
\gdef\SetFigFont#1#2#3{\begingroup
  \count@#1\relax \ifnum 25<\count@\count@25\fi
  \def\x{\endgroup\@setsize\SetFigFont{#2pt}}%
  \expandafter\x
    \csname \romannumeral\the\count@ pt\expandafter\endcsname
    \csname @\romannumeral\the\count@ pt\endcsname
  \csname #3\endcsname}%
\fi
\fi\endgroup
\begin{picture}(4821,1730)(526,-1046)
\put(526,-136){\makebox(0,0)[lb]{\smash{\SetFigFont{10}{12.0}{rm}{\color[rgb]{0,0,0}$v$}%
}}}
\put(908,-406){\makebox(0,0)[lb]{\smash{\SetFigFont{10}{12.0}{rm}{\color[rgb]{0,0,0}$1/2$}%
}}}
\put(1778,-375){\makebox(0,0)[lb]{\smash{\SetFigFont{10}{12.0}{rm}{\color[rgb]{0,0,0}$1/4$}%
}}}
\put(2408,-382){\makebox(0,0)[lb]{\smash{\SetFigFont{10}{12.0}{rm}{\color[rgb]{0,0,0}$3/4$}%
}}}
\put(1344,261){\makebox(0,0)[lb]{\smash{\SetFigFont{10}{12.0}{rm}{\color[rgb]{0,0,0}$1/2$}%
}}}
\put(2813,-90){\makebox(0,0)[lb]{\smash{\SetFigFont{10}{12.0}{rm}{\color[rgb]{0,0,0}$x$}%
}}}
\put(5347,-406){\makebox(0,0)[lb]{\smash{\SetFigFont{10}{12.0}{rm}{\color[rgb]{0,0,0}$w$}%
}}}
\put(3818,-436){\makebox(0,0)[lb]{\smash{\SetFigFont{10}{12.0}{rm}{\color[rgb]{0,0,0}$z$}%
}}}
\put(5213,179){\makebox(0,0)[lb]{\smash{\SetFigFont{10}{12.0}{rm}{\color[rgb]{0,0,0}$y$}%
}}}
\put(3152,-376){\makebox(0,0)[lb]{\smash{\SetFigFont{10}{12.0}{rm}{\color[rgb]{0,0,0}$3/8$}%
}}}
\put(3609,201){\makebox(0,0)[lb]{\smash{\SetFigFont{10}{12.0}{rm}{\color[rgb]{0,0,0}$3/8$}%
}}}
\put(4366,269){\makebox(0,0)[lb]{\smash{\SetFigFont{10}{12.0}{rm}{\color[rgb]{0,0,0}$3/8$}%
}}}
\put(2948,-849){\makebox(0,0)[lb]{\smash{\SetFigFont{10}{12.0}{rm}{\color[rgb]{0,0,0}$1/4$}%
}}}
\put(4651,-181){\makebox(0,0)[lb]{\smash{\SetFigFont{10}{12.0}{rm}{\color[rgb]{0,0,0}$5/8$}%
}}}
\put(1156,-997){\makebox(0,0)[lb]{\smash{\SetFigFont{10}{12.0}{rm}{\color[rgb]{0,0,0}$c=4$}%
}}}
\put(1254,-83){\makebox(0,0)[lb]{\smash{\SetFigFont{10}{12.0}{rm}{\color[rgb]{0,0,0}$w_1$}%
}}}
\put(2274,-105){\makebox(0,0)[lb]{\smash{\SetFigFont{10}{12.0}{rm}{\color[rgb]{0,0,0}$w_2$}%
}}}
\put(4553,-443){\makebox(0,0)[lb]{\smash{\SetFigFont{10}{12.0}{rm}{\color[rgb]{0,0,0}$w_4$}%
}}}
\end{picture}

%% file: EPS-A/redge2civ1b.pstex_t
\begin{picture}(0,0)%
\includegraphics{redge2civ1b.pstex}%
\end{picture}%
\setlength{\unitlength}{2763sp}%
\begingroup\makeatletter\ifx\SetFigFont\undefined
\def\x#1#2#3#4#5#6#7\relax{\def\x{#1#2#3#4#5#6}}%
\expandafter\x\fmtname xxxxxx\relax \def\y{splain}%
\ifx\x\y   
\gdef\SetFigFont#1#2#3{%
  \ifnum #1<17\tiny\else \ifnum #1<20\small\else
  \ifnum #1<24\normalsize\else \ifnum #1<29\large\else
  \ifnum #1<34\Large\else \ifnum #1<41\LARGE\else
     \huge\fi\fi\fi\fi\fi\fi
  \csname #3\endcsname}%
\else
\gdef\SetFigFont#1#2#3{\begingroup
  \count@#1\relax \ifnum 25<\count@\count@25\fi
  \def\x{\endgroup\@setsize\SetFigFont{#2pt}}%
  \expandafter\x
    \csname \romannumeral\the\count@ pt\expandafter\endcsname
    \csname @\romannumeral\the\count@ pt\endcsname
  \csname #3\endcsname}%
\fi
\fi\endgroup
\begin{picture}(4821,1687)(526,-1046)
\put(526,-136){\makebox(0,0)[lb]{\smash{\SetFigFont{8}{9.6}{rm}{\color[rgb]{0,0,0}$v$}%
}}}
\put(908,-406){\makebox(0,0)[lb]{\smash{\SetFigFont{8}{9.6}{rm}{\color[rgb]{0,0,0}$1/2$}%
}}}
\put(1778,-375){\makebox(0,0)[lb]{\smash{\SetFigFont{8}{9.6}{rm}{\color[rgb]{0,0,0}$1/4$}%
}}}
\put(2408,-382){\makebox(0,0)[lb]{\smash{\SetFigFont{8}{9.6}{rm}{\color[rgb]{0,0,0}$3/4$}%
}}}
\put(1344,261){\makebox(0,0)[lb]{\smash{\SetFigFont{8}{9.6}{rm}{\color[rgb]{0,0,0}$1/2$}%
}}}
\put(2813,-90){\makebox(0,0)[lb]{\smash{\SetFigFont{8}{9.6}{rm}{\color[rgb]{0,0,0}$x$}%
}}}
\put(5347,-406){\makebox(0,0)[lb]{\smash{\SetFigFont{8}{9.6}{rm}{\color[rgb]{0,0,0}$w$}%
}}}
\put(3818,-436){\makebox(0,0)[lb]{\smash{\SetFigFont{8}{9.6}{rm}{\color[rgb]{0,0,0}$z$}%
}}}
\put(5213,179){\makebox(0,0)[lb]{\smash{\SetFigFont{8}{9.6}{rm}{\color[rgb]{0,0,0}$y$}%
}}}
\put(3152,-376){\makebox(0,0)[lb]{\smash{\SetFigFont{8}{9.6}{rm}{\color[rgb]{0,0,0}$3/8$}%
}}}
\put(2948,-849){\makebox(0,0)[lb]{\smash{\SetFigFont{8}{9.6}{rm}{\color[rgb]{0,0,0}$1/4$}%
}}}
\put(4651,-181){\makebox(0,0)[lb]{\smash{\SetFigFont{8}{9.6}{rm}{\color[rgb]{0,0,0}$5/8$}%
}}}
\put(1156,-997){\makebox(0,0)[lb]{\smash{\SetFigFont{8}{9.6}{rm}{\color[rgb]{0,0,0}$c=4$}%
}}}
\put(5221,509){\makebox(0,0)[lb]{\smash{\SetFigFont{8}{9.6}{rm}{\color[rgb]{0,0,0}$t$}%
}}}
\put(3324,178){\makebox(0,0)[lb]{\smash{\SetFigFont{8}{9.6}{rm}{\color[rgb]{0,0,0}$3/8$}%
}}}
\put(4104,479){\makebox(0,0)[lb]{\smash{\SetFigFont{8}{9.6}{rm}{\color[rgb]{0,0,0}$3/16$}%
}}}
\put(4284,156){\makebox(0,0)[lb]{\smash{\SetFigFont{8}{9.6}{rm}{\color[rgb]{0,0,0}$3/16$}%
}}}
\put(1254,-83){\makebox(0,0)[lb]{\smash{\SetFigFont{8}{9.6}{rm}{\color[rgb]{0,0,0}$w_1$}%
}}}
\put(2274,-105){\makebox(0,0)[lb]{\smash{\SetFigFont{8}{9.6}{rm}{\color[rgb]{0,0,0}$w_2$}%
}}}
\put(4553,-443){\makebox(0,0)[lb]{\smash{\SetFigFont{8}{9.6}{rm}{\color[rgb]{0,0,0}$w_4$}%
}}}
\end{picture}

%% file: EPS-A/redge2civ2a.pstex_t
\begin{picture}(0,0)%
\includegraphics{redge2civ2a.pstex}%
\end{picture}%
\setlength{\unitlength}{2763sp}%
\begingroup\makeatletter\ifx\SetFigFont\undefined
\def\x#1#2#3#4#5#6#7\relax{\def\x{#1#2#3#4#5#6}}%
\expandafter\x\fmtname xxxxxx\relax \def\y{splain}%
\ifx\x\y   
\gdef\SetFigFont#1#2#3{%
  \ifnum #1<17\tiny\else \ifnum #1<20\small\else
  \ifnum #1<24\normalsize\else \ifnum #1<29\large\else
  \ifnum #1<34\Large\else \ifnum #1<41\LARGE\else
     \huge\fi\fi\fi\fi\fi\fi
  \csname #3\endcsname}%
\else
\gdef\SetFigFont#1#2#3{\begingroup
  \count@#1\relax \ifnum 25<\count@\count@25\fi
  \def\x{\endgroup\@setsize\SetFigFont{#2pt}}%
  \expandafter\x
    \csname \romannumeral\the\count@ pt\expandafter\endcsname
    \csname @\romannumeral\the\count@ pt\endcsname
  \csname #3\endcsname}%
\fi
\fi\endgroup
\begin{picture}(4410,1736)(526,-1046)
\put(526,-136){\makebox(0,0)[lb]{\smash{\SetFigFont{8}{9.6}{rm}{\color[rgb]{0,0,0}$v$}%
}}}
\put(908,-406){\makebox(0,0)[lb]{\smash{\SetFigFont{8}{9.6}{rm}{\color[rgb]{0,0,0}$1/2$}%
}}}
\put(1778,-375){\makebox(0,0)[lb]{\smash{\SetFigFont{8}{9.6}{rm}{\color[rgb]{0,0,0}$1/4$}%
}}}
\put(2408,-382){\makebox(0,0)[lb]{\smash{\SetFigFont{8}{9.6}{rm}{\color[rgb]{0,0,0}$3/4$}%
}}}
\put(1344,261){\makebox(0,0)[lb]{\smash{\SetFigFont{8}{9.6}{rm}{\color[rgb]{0,0,0}$1/2$}%
}}}
\put(2813,-90){\makebox(0,0)[lb]{\smash{\SetFigFont{8}{9.6}{rm}{\color[rgb]{0,0,0}$x$}%
}}}
\put(3818,-436){\makebox(0,0)[lb]{\smash{\SetFigFont{8}{9.6}{rm}{\color[rgb]{0,0,0}$z$}%
}}}
\put(3152,-376){\makebox(0,0)[lb]{\smash{\SetFigFont{8}{9.6}{rm}{\color[rgb]{0,0,0}$3/8$}%
}}}
\put(2948,-849){\makebox(0,0)[lb]{\smash{\SetFigFont{8}{9.6}{rm}{\color[rgb]{0,0,0}$1/4$}%
}}}
\put(3961,104){\makebox(0,0)[lb]{\smash{\SetFigFont{8}{9.6}{rm}{\color[rgb]{0,0,0}$3/8$}%
}}}
\put(3602,441){\makebox(0,0)[lb]{\smash{\SetFigFont{8}{9.6}{rm}{\color[rgb]{0,0,0}$3/8$}%
}}}
\put(4936,494){\makebox(0,0)[lb]{\smash{\SetFigFont{8}{9.6}{rm}{\color[rgb]{0,0,0}$y$}%
}}}
\put(4927, 97){\makebox(0,0)[lb]{\smash{\SetFigFont{8}{9.6}{rm}{\color[rgb]{0,0,0}$w$}%
}}}
\put(1156,-997){\makebox(0,0)[lb]{\smash{\SetFigFont{8}{9.6}{rm}{\color[rgb]{0,0,0}$c=27/8$}%
}}}
\put(1254,-83){\makebox(0,0)[lb]{\smash{\SetFigFont{8}{9.6}{rm}{\color[rgb]{0,0,0}$w_1$}%
}}}
\put(2274,-105){\makebox(0,0)[lb]{\smash{\SetFigFont{8}{9.6}{rm}{\color[rgb]{0,0,0}$w_2$}%
}}}
\put(4553,-443){\makebox(0,0)[lb]{\smash{\SetFigFont{8}{9.6}{rm}{\color[rgb]{0,0,0}$w_4$}%
}}}
\end{picture}

%% file: EPS-A/redge2civ2b.pstex_t
\begin{picture}(0,0)%
\includegraphics{redge2civ2b.pstex}%
\end{picture}%
\setlength{\unitlength}{2763sp}%
\begingroup\makeatletter\ifx\SetFigFont\undefined
\def\x#1#2#3#4#5#6#7\relax{\def\x{#1#2#3#4#5#6}}%
\expandafter\x\fmtname xxxxxx\relax \def\y{splain}%
\ifx\x\y   
\gdef\SetFigFont#1#2#3{%
  \ifnum #1<17\tiny\else \ifnum #1<20\small\else
  \ifnum #1<24\normalsize\else \ifnum #1<29\large\else
  \ifnum #1<34\Large\else \ifnum #1<41\LARGE\else
     \huge\fi\fi\fi\fi\fi\fi
  \csname #3\endcsname}%
\else
\gdef\SetFigFont#1#2#3{\begingroup
  \count@#1\relax \ifnum 25<\count@\count@25\fi
  \def\x{\endgroup\@setsize\SetFigFont{#2pt}}%
  \expandafter\x
    \csname \romannumeral\the\count@ pt\expandafter\endcsname
    \csname @\romannumeral\the\count@ pt\endcsname
  \csname #3\endcsname}%
\fi
\fi\endgroup
\begin{picture}(4702,2051)(526,-1046)
\put(526,-136){\makebox(0,0)[lb]{\smash{\SetFigFont{8}{9.6}{rm}{\color[rgb]{0,0,0}$v$}%
}}}
\put(908,-406){\makebox(0,0)[lb]{\smash{\SetFigFont{8}{9.6}{rm}{\color[rgb]{0,0,0}$1/2$}%
}}}
\put(1778,-375){\makebox(0,0)[lb]{\smash{\SetFigFont{8}{9.6}{rm}{\color[rgb]{0,0,0}$1/4$}%
}}}
\put(2408,-382){\makebox(0,0)[lb]{\smash{\SetFigFont{8}{9.6}{rm}{\color[rgb]{0,0,0}$3/4$}%
}}}
\put(1344,261){\makebox(0,0)[lb]{\smash{\SetFigFont{8}{9.6}{rm}{\color[rgb]{0,0,0}$1/2$}%
}}}
\put(2813,-90){\makebox(0,0)[lb]{\smash{\SetFigFont{8}{9.6}{rm}{\color[rgb]{0,0,0}$x$}%
}}}
\put(3818,-436){\makebox(0,0)[lb]{\smash{\SetFigFont{8}{9.6}{rm}{\color[rgb]{0,0,0}$z$}%
}}}
\put(3152,-376){\makebox(0,0)[lb]{\smash{\SetFigFont{8}{9.6}{rm}{\color[rgb]{0,0,0}$3/8$}%
}}}
\put(2948,-849){\makebox(0,0)[lb]{\smash{\SetFigFont{8}{9.6}{rm}{\color[rgb]{0,0,0}$1/4$}%
}}}
\put(1156,-997){\makebox(0,0)[lb]{\smash{\SetFigFont{8}{9.6}{rm}{\color[rgb]{0,0,0}$c=27/8$}%
}}}
\put(5228,846){\makebox(0,0)[lb]{\smash{\SetFigFont{8}{9.6}{rm}{\color[rgb]{0,0,0}$y$}%
}}}
\put(5219,104){\makebox(0,0)[lb]{\smash{\SetFigFont{8}{9.6}{rm}{\color[rgb]{0,0,0}$w$}%
}}}
\put(5228,479){\makebox(0,0)[lb]{\smash{\SetFigFont{8}{9.6}{rm}{\color[rgb]{0,0,0}$t$}%
}}}
\put(3730,718){\makebox(0,0)[lb]{\smash{\SetFigFont{8}{9.6}{rm}{\color[rgb]{0,0,0}$3/8$}%
}}}
\put(4403,156){\makebox(0,0)[lb]{\smash{\SetFigFont{8}{9.6}{rm}{\color[rgb]{0,0,0}$3/16$}%
}}}
\put(4149,449){\makebox(0,0)[lb]{\smash{\SetFigFont{8}{9.6}{rm}{\color[rgb]{0,0,0}$3/16$}%
}}}
\put(1254,-83){\makebox(0,0)[lb]{\smash{\SetFigFont{8}{9.6}{rm}{\color[rgb]{0,0,0}$w_1$}%
}}}
\put(2274,-105){\makebox(0,0)[lb]{\smash{\SetFigFont{8}{9.6}{rm}{\color[rgb]{0,0,0}$w_2$}%
}}}
\put(4553,-443){\makebox(0,0)[lb]{\smash{\SetFigFont{8}{9.6}{rm}{\color[rgb]{0,0,0}$w_4$}%
}}}
\end{picture}

%% file: EPS-A/redge2civ2c.pstex_t
\begin{picture}(0,0)%
\includegraphics{redge2civ2c.pstex}%
\end{picture}%
\setlength{\unitlength}{2763sp}%
\begingroup\makeatletter\ifx\SetFigFont\undefined
\def\x#1#2#3#4#5#6#7\relax{\def\x{#1#2#3#4#5#6}}%
\expandafter\x\fmtname xxxxxx\relax \def\y{splain}%
\ifx\x\y   
\gdef\SetFigFont#1#2#3{%
  \ifnum #1<17\tiny\else \ifnum #1<20\small\else
  \ifnum #1<24\normalsize\else \ifnum #1<29\large\else
  \ifnum #1<34\Large\else \ifnum #1<41\LARGE\else
     \huge\fi\fi\fi\fi\fi\fi
  \csname #3\endcsname}%
\else
\gdef\SetFigFont#1#2#3{\begingroup
  \count@#1\relax \ifnum 25<\count@\count@25\fi
  \def\x{\endgroup\@setsize\SetFigFont{#2pt}}%
  \expandafter\x
    \csname \romannumeral\the\count@ pt\expandafter\endcsname
    \csname @\romannumeral\the\count@ pt\endcsname
  \csname #3\endcsname}%
\fi
\fi\endgroup
\begin{picture}(4845,1681)(526,-1046)
\put(526,-136){\makebox(0,0)[lb]{\smash{\SetFigFont{8}{9.6}{rm}{\color[rgb]{0,0,0}$v$}%
}}}
\put(908,-406){\makebox(0,0)[lb]{\smash{\SetFigFont{8}{9.6}{rm}{\color[rgb]{0,0,0}$1/2$}%
}}}
\put(1778,-375){\makebox(0,0)[lb]{\smash{\SetFigFont{8}{9.6}{rm}{\color[rgb]{0,0,0}$1/4$}%
}}}
\put(2408,-382){\makebox(0,0)[lb]{\smash{\SetFigFont{8}{9.6}{rm}{\color[rgb]{0,0,0}$3/4$}%
}}}
\put(1344,261){\makebox(0,0)[lb]{\smash{\SetFigFont{8}{9.6}{rm}{\color[rgb]{0,0,0}$1/2$}%
}}}
\put(2813,-90){\makebox(0,0)[lb]{\smash{\SetFigFont{8}{9.6}{rm}{\color[rgb]{0,0,0}$x$}%
}}}
\put(3818,-436){\makebox(0,0)[lb]{\smash{\SetFigFont{8}{9.6}{rm}{\color[rgb]{0,0,0}$z$}%
}}}
\put(3152,-376){\makebox(0,0)[lb]{\smash{\SetFigFont{8}{9.6}{rm}{\color[rgb]{0,0,0}$3/8$}%
}}}
\put(2948,-849){\makebox(0,0)[lb]{\smash{\SetFigFont{8}{9.6}{rm}{\color[rgb]{0,0,0}$1/4$}%
}}}
\put(5219,104){\makebox(0,0)[lb]{\smash{\SetFigFont{8}{9.6}{rm}{\color[rgb]{0,0,0}$w$}%
}}}
\put(5371,-399){\makebox(0,0)[lb]{\smash{\SetFigFont{8}{9.6}{rm}{\color[rgb]{0,0,0}$t$}%
}}}
\put(5228,479){\makebox(0,0)[lb]{\smash{\SetFigFont{8}{9.6}{rm}{\color[rgb]{0,0,0}$y$}%
}}}
\put(4592,-173){\makebox(0,0)[lb]{\smash{\SetFigFont{8}{9.6}{rm}{\color[rgb]{0,0,0}$7/16$}%
}}}
\put(3908,-182){\makebox(0,0)[lb]{\smash{\SetFigFont{8}{9.6}{rm}{\color[rgb]{0,0,0}$3/16$}%
}}}
\put(4110,216){\makebox(0,0)[lb]{\smash{\SetFigFont{8}{9.6}{rm}{\color[rgb]{0,0,0}$3/16$}%
}}}
\put(3633,403){\makebox(0,0)[lb]{\smash{\SetFigFont{8}{9.6}{rm}{\color[rgb]{0,0,0}$3/8$}%
}}}
\put(1156,-997){\makebox(0,0)[lb]{\smash{\SetFigFont{8}{9.6}{rm}{\color[rgb]{0,0,0}$c=61/16$}%
}}}
\put(1254,-83){\makebox(0,0)[lb]{\smash{\SetFigFont{8}{9.6}{rm}{\color[rgb]{0,0,0}$w_1$}%
}}}
\put(2274,-105){\makebox(0,0)[lb]{\smash{\SetFigFont{8}{9.6}{rm}{\color[rgb]{0,0,0}$w_2$}%
}}}
\put(4553,-443){\makebox(0,0)[lb]{\smash{\SetFigFont{8}{9.6}{rm}{\color[rgb]{0,0,0}$w_4$}%
}}}
\end{picture}

%% file: other.tex

\section{Other Pivot Rules}
\label{sec:other}


\subsection{Bland's Rule}
\label{subsec:leastInd}
\label{subsec:bland}

For Bland's \emph{least index} pivot rule~\cite{bland77:_new} the
facets (inequalities) are numbered.  At every non-minimal vertex the
rule then dictates to choose the edge that leaves the facet with the
smallest number.  
(A special feature of Bland's rule is that it does not admit
cycling even on degenerate programs/\allowbreak non-simple polytopes,
when our geometric description of the rule is, however, not applicable.)

\begin{proposition}
  \label{prop:bland}
  The linearity coefficient of Bland's rule is~$2$. 
\end{proposition}

\begin{proof}
  Figure \ref{fig:li} illustrates a family of $3$-dimensional LPs
  on which Bland's rule, started at $v_{\text{start}} = v_{2n-6}$,
  visits all but one of the vertices. (As we have already noted, the directed graph in the figure is
  readily verified to satisfy the conditions of Theorem~\ref{thm:MK}.)
  Specifically,
  choose an initial numbering of the facets where the largest index is
  assigned to facet $f$. When starting at the vertex $v_\st = v_{2n-6}$
  the   simplex algorithm with Bland's rule visits the $2n-5$
  vertices $v_{2n-6}, \ldots, v_0$.
\end{proof}


\subsection{Dantzig's Rule}
\label{subsec:largestCoeff}
\label{subsect:dantzig}

\emph{Dantzig's rule} is the original rule proposed by Dantzig
when he invented the simplex algorithm. 
In his setting of
a maximization problem formulated in the language of simplex tableaus,
the rule requires to pivot into the basis the variable that has the
largest reduced cost coefficient (if no variable has positive
reduced cost, the current tableau is optimal).

By suitably scaling the inequalities of the LP, Dantzig's rule follows
the same path as Bland's rule; see Amenta \& Ziegler
\cite[Observation 2.6]{AZ96}. Thus Dantzig's rule cannot be
faster than Bland's rule, and Proposition \ref{prop:bland} thus implies:

\begin{proposition}
  \label{prop:dantzig}
  The linearity coefficient of Dantzig's rule is~$2$. 
\end{proposition}


\subsection{Greatest Decrease Rule}
\label{subsec:greatestDecr}

The \emph{greatest decrease} rule moves from any non-optimal
vertex to the neighbor with the smallest objective function value. 
We assume that the objective function is generic, so the vertex is
unique. However, the greatest decrease rule may compare
non-adjacent neighbors, so the information given by
the directed graph is not sufficient to implement it; we rather need explicit objective function values.

\begin{proposition}
The linearity coefficient of the greatest decrease rule is~$\frac32$.
\end{proposition}

\begin{proof}
First we show that $\Lambda(\GD)\ge\frac32$.
Figure \ref{fig:gd} indicates a family of $3$-dimensional LPs. 
By Theorem \ref{thm:MK}, there is a realization of these LPs 
with the objective function 
linear 
ordering on the vertices given 
by the left-to-right ordering in our figure. 
Started at $v_{\st} = v_{2n-6}$, the greatest decrease
rule visits all $1$-vertices, the global sink, and half of the
$2$-vertices. Thus it needs $\frac32 (n-3)$ pivot 
steps to reach $v_{\min} = v_0$.

  \begin{figure}[ht]
    \begin{center}
      \input{EPS/gd.pstex_t}
      \caption{ \label{fig:gd}
        Lower bound for the greatest decrease rule. 
        All edges are oriented from left to right.
      }
    \end{center}
  \end{figure}
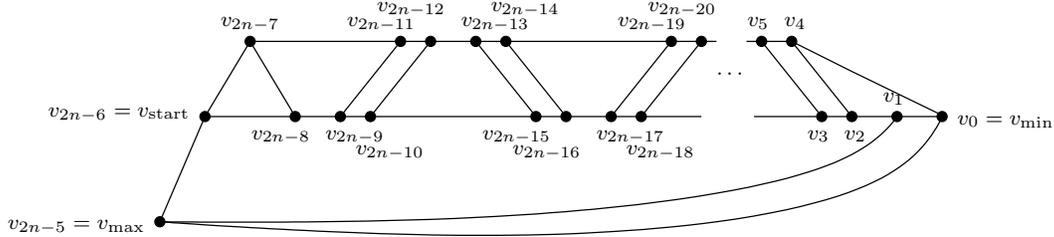

\noindent
For the proof of $\Lambda(\GD)\le\frac32$, we consider
an arbitrary instance with $n$, $P$, $\varphi$, and $v_{\st}$ as above.
Denote by $n_1$ and $n_2$ the number of visited $1$- and $2$-vertices, respectively. 
Thus there are $n-3-n_1$ and
$n-3-n_2$ unvisited $1$-vertices and $2$-vertices, respectively.
For every visited $2$-vertex $v$ only one of the two direct
successors $v'$ and $v''$ is visited. Assuming that $\varphi(v')>\varphi(v'')$,
the greatest decrease rule will proceed directly from $v$ to~$v''$ and thus
skip $v'$, whose objective function value satisfies
$\varphi(v)>\varphi(v')>\varphi(v'')$. Thus there is an unvisited vertex
uniquely associated with every visited $2$-vertex. Thus
$n_2\le 2n-6-n_1-n_2$,
which is equivalent to $n_1 + 2n_2\le 2n-6$. We get
\[n_1+n_2\ =\ \tfrac 12 n_1 +\tfrac 12 (n_1 + 2n_2)
\ \le\ \tfrac 12 (n-3) + \tfrac 12 (2n-6)\ \le\ \tfrac 32 (n-3).
\]
This yields $\Lambda(\GD)\le\frac32$ and completes the proof.
\end{proof}


\subsection{Steepest Decrease Rule}
\label{subsec:steepestDecr}
At any non-minimal vertex $v$ the \emph{steepest decrease} pivot rule
moves to the neighbor $w$ with $vw$ 
being the steepest decreasing
edge, that is, such that $\frac {\langle c,w-v\rangle} {\lVert w-v
  \rVert \, \lVert c \rVert}$ is minimal 
(where $\langle c,x\rangle$ is the objective function).

\begin{proposition} \label{prop:sd}
The linearity coefficient of the steepest decrease rule is~$2$.
\end{proposition}

\begin{figure}[ht]
    \begin{center}
      \input{EPS/sd-new.pstex_t}
      \caption{ \label{fig:sd}
        Lower bounds for the steepest decrease and shadow vertex rules.
        Planar projection of the polytope: The objective function is
        $x_1$; it directs all edges from left to~right.}
    \end{center}
  \end{figure}
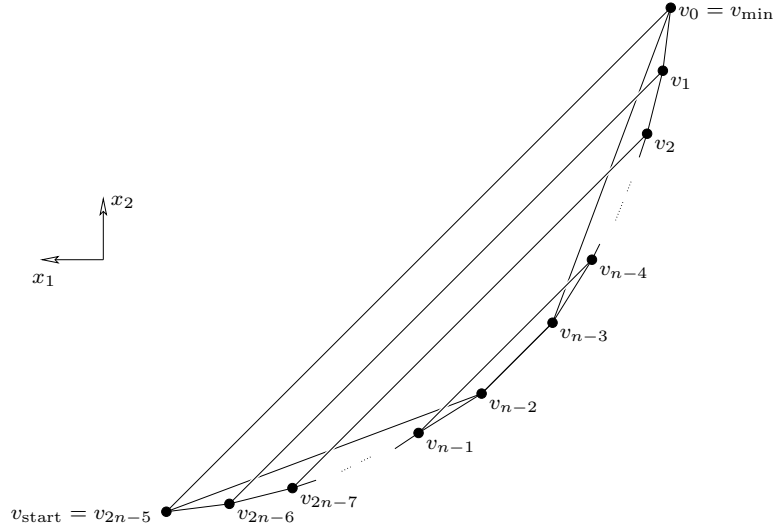

\begin{proof}
Figure \ref{fig:sd} depicts a planar projection
onto the $(x_1,x_2)$-plane of an LP that
is easily constructed either ``by hand'' or as a deformed
product (see Amenta \& Ziegler \cite{AZ96}).
If the polytope is scaled to be very flat in the $x_3$-direction,
then steepest decrease tells the simplex algorithm to use the
edge that in the projection has the smallest slope (in absolute value).
Thus starting at $v_{\st} = v_{2n-5}$, the steepest decrease
rule visits \emph{all} the vertices.
\end{proof}

\subsection{Shadow Vertex Rule}
\label{subsec:shadow}

The \emph{shadow vertex} pivot rule chooses a sequence of
edges that lie on the boundary of the $2$-dimensional projection
of the polytope given by $x\mapsto (\langle c,x\rangle,\langle d,x\rangle)$, where
$\langle c,x\rangle$ is the given objective function, and 
$\langle d,x\rangle$ is an objective function that is constructed
to be optimal at the starting vertex $v_\st$.
The vertices that are visited on the path from $v_\st$ to
$v_{\min}$ are then optimal for objective functions that interpolate
between $\langle d,x\rangle$ and $\langle c,x\rangle$.
(This pivot rule is known to be polynomial on ``random linear programs''
in specific models; cf. Borgwardt \cite{Borg}, Ziegler \cite{Z79},
and Spielman \& Teng \cite{spielman02:_smoot_analy_algor}.)

\begin{proposition}
The linearity coefficient of the shadow vertex rule is~$2$.
\end{proposition}

\begin{proof}
We reuse the linear programs of Proposition \ref{prop:sd}/Figure \ref{fig:sd}.
Here $v_{2n-5}=v_{\max}$ is optimal for the
starting objective function $\langle d,x\rangle=x_2$, while $v_0$ is optimal
for $\langle c,x\rangle=x_1$. On the way from $v_{2n-5}$ to $v_{\min}=v_0$ the
shadow vertex rule visits \emph{all} the vertices.
\end{proof}


\subsection{Random Facet}
\label{subsec:RF}

The 
\emph{random facet} pivot rule, due to Kalai \cite[p.~228]{Ka97}, is as
follows:
\begin{itemize}
\item[(RF)] At any non-optimal vertex $v$ choose one facet $f$
  containing $v$ uniformly at random
  and solve the problem restricted to $f$ by applying (RF) recursively.\\
  The recursion will eventually restrict to a one-dimensional
  subproblem (that is, an edge), which is solved by following the
  edge.
\end{itemize}
The one-dimensional base case singled out here is only
  implicit in Kalai's work. This is probably the reason why
there are different versions of this rule in the literature which 
unfortunately were
not distinguished. They all differ in the way how 
1-vertices are treated. Since the (unique) out-edge of a $1$-vertex is
always taken with probability one (regardless of which facets we
restrict to) we could use the following
alternative formulations of the random facet rule:
\begin{itemize}
\item[(RF1)]
At each non-optimal vertex $v$ follow the (unique) outgoing edge if
$v$ is a 1-vertex. Otherwise choose one facet $f$ uniformly at random
containing $v$ and 
solve the problem restricted to $f$ by applying (RF1)
recursively. 
\item[(RF2)]
At any non-optimal vertex $v$ choose one facet $f$ containing $v$
uniformly at random 
and solve the problem restricted to $f$ by applying 
(RF2) recursively. 
The minimal vertex $\operatorname{opt}(f)$ of $f$ is a 1-vertex and we follow
the (unique) outgoing edge of the vertex $\operatorname{opt}(f)$.
\end{itemize}
The variant (RF1) appears in G\"artner, Henk \& Ziegler \cite[p.~350]{GHZ98},
while the version (RF2) is from G\"artner \cite{Ga98},
who, however, formulated this variant of the
random facet rule for combinatorial cubes, where
the formulations above are equivalent.

Note that (RF) uses randomness at every vertex,
and (RF1) would follow a path of $1$-vertices
deterministically, while (RF2) takes at most one deterministic
step in a row. This results in distinct pivot rules,
with different worst case examples. 

\begin{proposition}\label{prop:lf}
For each version (RF), (RF1), and (RF2) of the
random facet rule the linearity coefficient is~$2$.
\end{proposition}

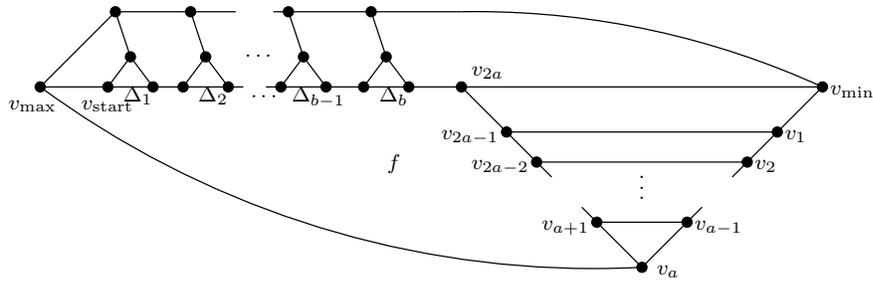
\begin{figure}
    \begin{center}
      \input{EPS/le-re.pstex_t}
      \caption{ \label{fig:le-re}
        Lower bounds for the (RF1) and (RF2) variants of the random
        facet rule, and for the least entered rule with random edge as
        the tie-breaking rule.  All edges are oriented from left to
        right.  }
    \end{center}
  \end{figure}

\begin{proof}
Figure \ref{fig:rf} depicts a family of LPs with 
$2n-4=2a+2b+2$ vertices and $n = a + b + 3$
facets. For each of the $b$ $1$-vertices 
$v_\st=v_{2n-7}, v_{2n-9}, \dots, v_{2a+1}$, 
the probability of leaving it via choosing facet~$f$ is   
$\tfrac 12$. After choosing facet~$f$, (RF) ``sticks'' to facet $f$ until $v_a$ is reached. 

Choosing $a = k^2$ and $b= k$ we obtain a family of LPs with 
$n =k^2 + k + 3$ facets. Then (RF) sticks to facet $f$ with
probability $p \ge 1-(\frac 12)^k$.
Thus the expected number of visited vertices is at least 
\[
\big(1-(\tfrac 12)^k\big)\,(2a+b)\ \ge\ 2 k^2 - \frac {2 k^2}{2^k}.
\]
Since there are $n=k^2+k+3$ facets, the linearity coefficient is $2$.

  \begin{figure}[ht]
    \begin{center}
      \input{EPS/rf.pstex_t}
      \caption{ \label{fig:rf}
        Lower bound for the random facet rule (RF).
        All edges are oriented from left to right.
      }
    \end{center}
  \end{figure}
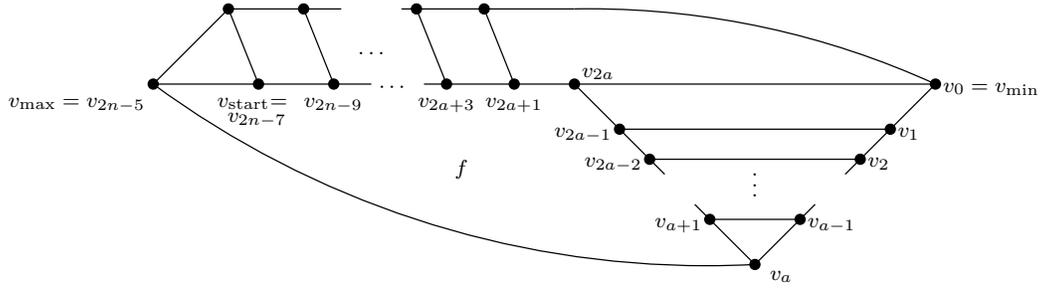
  
  The version (RF1) of the random facet rule follows the path of $1$-vertices 
  $v_\st=v_{2n-7}, v_{2n-9}, \dots, v_{2a+1}$ deterministically. 
  We can cut off each of these vertices. 
  This yields the graphs depicted in Figure
  \ref{fig:le-re}. At each source of the new facets
  $\Delta_1,\ldots,\Delta_b$, the facet $f$ is chosen with
  probability $\tfrac 13$. If any of the other two facets is chosen, we end up at the sink vertex of the respective facet $\Delta_i$. Thus the linearity coefficient remains~$2$,
  only the rate of convergence decreases. 
  The same works for (RF2) as well.
\end{proof}

\subsection{Least Entered Rule}
\label{subsec:leastEnt}
At any non-optimal vertex, the \emph{least entered}
pivot rule chooses the decreasing edge that
leaves the facet that has been left least often in the previous
moves. In case of ties a tie-breaking rule is used to determine the
decreasing edge to be taken. Any other pivot rule 
can be used as a tie-breaking rule.

The least entered rule was first formulated by Norman Zadeh around 1980
(see \cite{KlKl1} and \cite{Z79}). 
It has still not been determined whether Zadeh's rule is
polynomial if the dimension is part of the input. 
Zadeh has offered \$1000 for solving this problem.

\begin{proposition}
The linearity coefficient of the least entered rule with greatest
decrease as tie-breaking rule is~$2$.
\end{proposition}

  \begin{figure}[ht]
    \begin{center}
      \input{EPS/le-gd.pstex_t}
      \caption{\label{fig:le-gd}
        Lower bound for the least entered rule
        with greatest decrease as the tie-breaking rule.
        All edges are oriented from left to right.
      }
    \end{center}
  \end{figure}
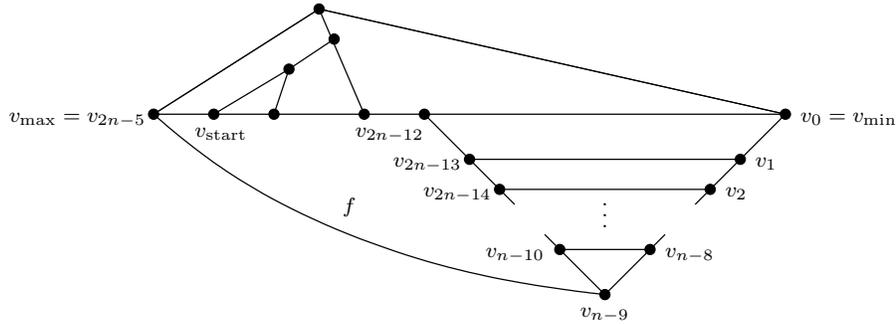

\begin{proof}
  Figure \ref{fig:le-gd} describes a family of $3$-dimensional LPs,
  where the left-to-right ordering of the vertices suggested
  by the figure can be realized, according to the Mihalisin--Klee
  Theorem \ref{thm:MK}.
  Starting at $v_{\st}=v_{2n-6}$, the greatest decrease rule decides to leave
  the facet $f$. Following two $1$-vertices the facet $f$ is entered
  again. All upcoming facets have not been visited before. Thus the
  least entered rule ``sticks'' to the facet $f$ and $2n-7$ vertices
  (that is, all but $3$ vertices) are visited.
\end{proof}

\begin{proposition} \label{prop:le-re}
The linearity coefficient of the least entered rule with 
random edge as the tie-breaking rule is~$2$.
\end{proposition}


\begin{proof}
  Figure \ref{fig:le-re} describes LPs with $2n-4=2a+4b+2$ vertices and 
  $n = a + 2b+3$ facets.  At the sources of the facets
  $\Delta_i$ the random edge rule leaves the facet $f$ with
  probability $\tfrac 12$. As soon as $f$ is left once, it will be
  revisited and the least entered rule will ``stick'' to the facet~$f$. 
  (Thus the only way not to ``stick'' to $f$ is that the random edge rule
  chooses to continue along $f$ until it reaches
  the vertex $v_{2a}$.)
  When the least entered rule ``sticks'' to the facet $f$ all of the $2a$ vertices
  $v_{2a-1}, v_{2a-2},\dots,v_1,v_0$ are visited.

Now the analysis is
exactly the same as in the proof of Proposition \ref{prop:lf}. Thus
choosing $a=k^2$ and $b=k$ yields that the linearity
coefficient is $2$.
\end{proof}







%% file: EPS/gd.pstex_t
\begin{picture}(0,0)%
\includegraphics{gd.pstex}%
\end{picture}%
\setlength{\unitlength}{4144sp}%
\begingroup\makeatletter\ifx\SetFigFont\undefined
\def\x#1#2#3#4#5#6#7\relax{\def\x{#1#2#3#4#5#6}}%
\expandafter\x\fmtname xxxxxx\relax \def\y{splain}%
\ifx\x\y   
\gdef\SetFigFont#1#2#3{%
  \ifnum #1<17\tiny\else \ifnum #1<20\small\else
  \ifnum #1<24\normalsize\else \ifnum #1<29\large\else
  \ifnum #1<34\Large\else \ifnum #1<41\LARGE\else
     \huge\fi\fi\fi\fi\fi\fi
  \csname #3\endcsname}%
\else
\gdef\SetFigFont#1#2#3{\begingroup
  \count@#1\relax \ifnum 25<\count@\count@25\fi
  \def\x{\endgroup\@setsize\SetFigFont{#2pt}}%
  \expandafter\x
    \csname \romannumeral\the\count@ pt\expandafter\endcsname
    \csname @\romannumeral\the\count@ pt\endcsname
  \csname #3\endcsname}%
\fi
\fi\endgroup
\begin{picture}(4860,1437)(361,-694)
\put(3871,254){\makebox(0,0)[b]{\smash{\SetFigFont{8}{9.6}{rm}{\color[rgb]{0,0,0}$\cdots$}%
}}}
\put(2161,659){\makebox(0,0)[rb]{\smash{\SetFigFont{8}{9.6}{rm}{\color[rgb]{0,0,0}$v_{2n-12}$}%
}}}
\put(991,569){\makebox(0,0)[b]{\smash{\SetFigFont{8}{9.6}{rm}{\color[rgb]{0,0,0}$v_{2n-7}$}%
}}}
\put(1981,569){\makebox(0,0)[rb]{\smash{\SetFigFont{8}{9.6}{rm}{\color[rgb]{0,0,0}$v_{2n-11}$}%
}}}
\put(2251,569){\makebox(0,0)[lb]{\smash{\SetFigFont{8}{9.6}{rm}{\color[rgb]{0,0,0}$v_{2n-13}$}%
}}}
\put(2431,659){\makebox(0,0)[lb]{\smash{\SetFigFont{8}{9.6}{rm}{\color[rgb]{0,0,0}$v_{2n-14}$}%
}}}
\put(1351,-106){\makebox(0,0)[rb]{\smash{\SetFigFont{8}{9.6}{rm}{\color[rgb]{0,0,0}$v_{2n-8}$}%
}}}
\put(1441,-106){\makebox(0,0)[lb]{\smash{\SetFigFont{8}{9.6}{rm}{\color[rgb]{0,0,0}$v_{2n-9}$}%
}}}
\put(1621,-196){\makebox(0,0)[lb]{\smash{\SetFigFont{8}{9.6}{rm}{\color[rgb]{0,0,0}$v_{2n-10}$}%
}}}
\put(3061,-106){\makebox(0,0)[lb]{\smash{\SetFigFont{8}{9.6}{rm}{\color[rgb]{0,0,0}$v_{2n-17}$}%
}}}
\put(2791,-106){\makebox(0,0)[rb]{\smash{\SetFigFont{8}{9.6}{rm}{\color[rgb]{0,0,0}$v_{2n-15}$}%
}}}
\put(2971,-196){\makebox(0,0)[rb]{\smash{\SetFigFont{8}{9.6}{rm}{\color[rgb]{0,0,0}$v_{2n-16}$}%
}}}
\put(3241,-196){\makebox(0,0)[lb]{\smash{\SetFigFont{8}{9.6}{rm}{\color[rgb]{0,0,0}$v_{2n-18}$}%
}}}
\put(3601,569){\makebox(0,0)[rb]{\smash{\SetFigFont{8}{9.6}{rm}{\color[rgb]{0,0,0}$v_{2n-19}$}%
}}}
\put(3781,659){\makebox(0,0)[rb]{\smash{\SetFigFont{8}{9.6}{rm}{\color[rgb]{0,0,0}$v_{2n-20}$}%
}}}
\put(4906,119){\makebox(0,0)[rb]{\smash{\SetFigFont{8}{9.6}{rm}{\color[rgb]{0,0,0}$v_1$}%
}}}
\put(4096,569){\makebox(0,0)[rb]{\smash{\SetFigFont{8}{9.6}{rm}{\color[rgb]{0,0,0}$v_5$}%
}}}
\put(4186,569){\makebox(0,0)[lb]{\smash{\SetFigFont{8}{9.6}{rm}{\color[rgb]{0,0,0}$v_4$}%
}}}
\put(4456,-106){\makebox(0,0)[rb]{\smash{\SetFigFont{8}{9.6}{rm}{\color[rgb]{0,0,0}$v_3$}%
}}}
\put(4546,-106){\makebox(0,0)[lb]{\smash{\SetFigFont{8}{9.6}{rm}{\color[rgb]{0,0,0}$v_2$}%
}}}
\put(361,-646){\makebox(0,0)[rb]{\smash{\SetFigFont{8}{9.6}{rm}{\color[rgb]{0,0,0}$v_{2n-5} = v_{\max}$}%
}}}
\put(5221,-16){\makebox(0,0)[lb]{\smash{\SetFigFont{8}{9.6}{rm}{\color[rgb]{0,0,0}$v_0=v_{\min}$}%
}}}
\put(631, 29){\makebox(0,0)[rb]{\smash{\SetFigFont{8}{9.6}{rm}{\color[rgb]{0,0,0}$v_{2n-6}= v_\st$}%
}}}
\end{picture}

%% file: EPS/sd-new.pstex_t
\begin{picture}(0,0)%
\includegraphics{sd-new.pstex}%
\end{picture}%
\setlength{\unitlength}{1302sp}%
\begingroup\makeatletter\ifx\SetFigFont\undefined
\def\x#1#2#3#4#5#6#7\relax{\def\x{#1#2#3#4#5#6}}%
\expandafter\x\fmtname xxxxxx\relax \def\y{splain}%
\ifx\x\y   
\gdef\SetFigFont#1#2#3{%
  \ifnum #1<17\tiny\else \ifnum #1<20\small\else
  \ifnum #1<24\normalsize\else \ifnum #1<29\large\else
  \ifnum #1<34\Large\else \ifnum #1<41\LARGE\else
     \huge\fi\fi\fi\fi\fi\fi
  \csname #3\endcsname}%
\else
\gdef\SetFigFont#1#2#3{\begingroup
  \count@#1\relax \ifnum 25<\count@\count@25\fi
  \def\x{\endgroup\@setsize\SetFigFont{#2pt}}%
  \expandafter\x
    \csname \romannumeral\the\count@ pt\expandafter\endcsname
    \csname @\romannumeral\the\count@ pt\endcsname
  \csname #3\endcsname}%
\fi
\fi\endgroup
\begin{picture}(12300,9996)(-149,-10219)
\put(3751,-10111){\makebox(0,0)[lb]{\smash{\SetFigFont{8}{9.6}{rm}{\color[rgb]{0,0,0}$v_{2n-6}$}%
}}}
\put(4951,-9811){\makebox(0,0)[lb]{\smash{\SetFigFont{8}{9.6}{rm}{\color[rgb]{0,0,0}$v_{2n-7}$}%
}}}
\put(9901,-6661){\makebox(0,0)[lb]{\smash{\SetFigFont{8}{9.6}{rm}{\color[rgb]{0,0,0}$v_{n-3}$}%
}}}
\put(8551,-8011){\makebox(0,0)[lb]{\smash{\SetFigFont{8}{9.6}{rm}{\color[rgb]{0,0,0}$v_{n-2}$}%
}}}
\put(7351,-8761){\makebox(0,0)[lb]{\smash{\SetFigFont{8}{9.6}{rm}{\color[rgb]{0,0,0}$v_{n-1}$}%
}}}
\put(10651,-5461){\makebox(0,0)[lb]{\smash{\SetFigFont{8}{9.6}{rm}{\color[rgb]{0,0,0}$v_{n-4}$}%
}}}
\put(11701,-3061){\makebox(0,0)[lb]{\smash{\SetFigFont{8}{9.6}{rm}{\color[rgb]{0,0,0}$v_2$}%
}}}
\put(12001,-1861){\makebox(0,0)[lb]{\smash{\SetFigFont{8}{9.6}{rm}{\color[rgb]{0,0,0}$v_1$}%
}}}
\put(12151,-511){\makebox(0,0)[lb]{\smash{\SetFigFont{8}{9.6}{rm}{\color[rgb]{0,0,0}$v_0=v_{\min}$}%
}}}
\put(-149,-5611){\makebox(0,0)[lb]{\smash{\SetFigFont{8}{9.6}{rm}{\color[rgb]{0,0,0}$x_1$}%
}}}
\put(1351,-4111){\makebox(0,0)[lb]{\smash{\SetFigFont{8}{9.6}{rm}{\color[rgb]{0,0,0}$x_2$}%
}}}
\put(2146,-10111){\makebox(0,0)[rb]{\smash{\SetFigFont{8}{9.6}{rm}{\color[rgb]{0,0,0}$v_\st=v_{2n-5}$}%
}}}
\end{picture}

%% file: EPS/rf.pstex_t
\begin{picture}(0,0)%
\includegraphics{rf.pstex}%
\end{picture}%
\setlength{\unitlength}{4144sp}%
\begingroup\makeatletter\ifx\SetFigFont\undefined%
\gdef\SetFigFont#1#2#3#4#5{%
  \reset@font\fontsize{#1}{#2pt}%
  \fontfamily{#3}\fontseries{#4}\fontshape{#5}%
  \selectfont}%
\fi\endgroup%
\begin{picture}(4807,1687)(496,-902)
\put(926,164){\makebox(0,0)[lb]{\smash{{\SetFigFont{8}{9.6}{\updefault}{\color[rgb]{0,0,0}$v_\st{=}$}%
}}}}
\put(986, 74){\makebox(0,0)[lb]{\smash{{\SetFigFont{8}{9.6}{\updefault}{\color[rgb]{0,0,0}$v_{2n-7}$}%
}}}}
\put(4141,-376){\makebox(0,0)[b]{\smash{{\SetFigFont{8}{9.6}{\updefault}{\color[rgb]{0,0,0}$\vdots$}%
}}}}
\put(5266,254){\makebox(0,0)[lb]{\smash{{\SetFigFont{8}{9.6}{\updefault}{\color[rgb]{0,0,0}$v_0=v_{\min}$}%
}}}}
\put(4816,-196){\makebox(0,0)[lb]{\smash{{\SetFigFont{8}{9.6}{\updefault}{\color[rgb]{0,0,0}$v_2$}%
}}}}
\put(4456,-556){\makebox(0,0)[lb]{\smash{{\SetFigFont{8}{9.6}{\updefault}{\color[rgb]{0,0,0}$v_{a-1}$}%
}}}}
\put(4231,-871){\makebox(0,0)[lb]{\smash{{\SetFigFont{8}{9.6}{\updefault}{\color[rgb]{0,0,0}$v_a$}%
}}}}
\put(3826,-556){\makebox(0,0)[rb]{\smash{{\SetFigFont{8}{9.6}{\updefault}{\color[rgb]{0,0,0}$v_{a+1}$}%
}}}}
\put(3466,-196){\makebox(0,0)[rb]{\smash{{\SetFigFont{8}{9.6}{\updefault}{\color[rgb]{0,0,0}$v_{2a-2}$}%
}}}}
\put(3286,-16){\makebox(0,0)[rb]{\smash{{\SetFigFont{8}{9.6}{\updefault}{\color[rgb]{0,0,0}$v_{2a-1}$}%
}}}}
\put(496,164){\makebox(0,0)[rb]{\smash{{\SetFigFont{8}{9.6}{\updefault}{\color[rgb]{0,0,0}$v_{\max}=v_{2n-5}$}%
}}}}
\put(1846,479){\makebox(0,0)[b]{\smash{{\SetFigFont{8}{9.6}{\updefault}{\color[rgb]{0,0,0}$\ldots$}%
}}}}
\put(1981,254){\makebox(0,0)[b]{\smash{{\SetFigFont{8}{9.6}{\updefault}{\color[rgb]{0,0,0}$\cdots$}%
}}}}
\put(2296,164){\makebox(0,0)[b]{\smash{{\SetFigFont{8}{9.6}{\updefault}{\color[rgb]{0,0,0}$v_{2a+3}$}%
}}}}
\put(4996,-16){\makebox(0,0)[lb]{\smash{{\SetFigFont{8}{9.6}{\updefault}{\color[rgb]{0,0,0}$v_1$}%
}}}}
\put(2701,164){\makebox(0,0)[b]{\smash{{\SetFigFont{8}{9.6}{\updefault}{\color[rgb]{0,0,0}$v_{2a+1}$}%
}}}}
\put(3106,344){\makebox(0,0)[lb]{\smash{{\SetFigFont{8}{9.6}{\updefault}{\color[rgb]{0,0,0}$v_{2a}$}%
}}}}
\put(2386,-241){\makebox(0,0)[b]{\smash{{\SetFigFont{8}{9.6}{\updefault}{\color[rgb]{0,0,0}$f$}%
}}}}
\put(1621,164){\makebox(0,0)[b]{\smash{{\SetFigFont{8}{9.6}{\updefault}{\color[rgb]{0,0,0}$v_{2n-9}$}%
}}}}
\end{picture}%

%% file: EPS/le-re.pstex_t
\begin{picture}(0,0)%
\includegraphics{le-re.pstex}%
\end{picture}%
\setlength{\unitlength}{4144sp}%
\begingroup\makeatletter\ifx\SetFigFont\undefined%
\gdef\SetFigFont#1#2#3#4#5{%
  \reset@font\fontsize{#1}{#2pt}%
  \fontfamily{#3}\fontseries{#4}\fontshape{#5}%
  \selectfont}%
\fi\endgroup%
\begin{picture}(4807,1642)(496,-857)
\put(5266,254){\makebox(0,0)[lb]{\smash{{\SetFigFont{8}{9.6}{\updefault}{\color[rgb]{0,0,0}$v_{\text{min}}$}%
}}}}
\put(4141,-376){\makebox(0,0)[b]{\smash{{\SetFigFont{8}{9.6}{\updefault}{\color[rgb]{0,0,0}$\vdots$}%
}}}}
\put(2656,-196){\makebox(0,0)[b]{\smash{{\SetFigFont{8}{9.6}{\updefault}{\color[rgb]{0,0,0}$f$}%
}}}}
\put(496,164){\makebox(0,0)[b]{\smash{{\SetFigFont{8}{9.6}{\updefault}{\color[rgb]{0,0,0}$v_{\text{max}}$}%
}}}}
\put(946,164){\makebox(0,0)[b]{\smash{{\SetFigFont{8}{9.6}{\updefault}{\color[rgb]{0,0,0}$v_{\text{start}}$}%
}}}}
\put(2656,209){\makebox(0,0)[b]{\smash{{\SetFigFont{8}{9.6}{\updefault}{\color[rgb]{0,0,0}$\Delta_b$}%
}}}}
\put(4996,-16){\makebox(0,0)[lb]{\smash{{\SetFigFont{8}{9.6}{\updefault}{\color[rgb]{0,0,0}$v_1$}%
}}}}
\put(4816,-196){\makebox(0,0)[lb]{\smash{{\SetFigFont{8}{9.6}{\updefault}{\color[rgb]{0,0,0}$v_2$}%
}}}}
\put(4231,-826){\makebox(0,0)[lb]{\smash{{\SetFigFont{8}{9.6}{\updefault}{\color[rgb]{0,0,0}$v_a$}%
}}}}
\put(4456,-556){\makebox(0,0)[lb]{\smash{{\SetFigFont{8}{9.6}{\updefault}{\color[rgb]{0,0,0}$v_{a-1}$}%
}}}}
\put(3826,-556){\makebox(0,0)[rb]{\smash{{\SetFigFont{8}{9.6}{\updefault}{\color[rgb]{0,0,0}$v_{a+1}$}%
}}}}
\put(3466,-196){\makebox(0,0)[rb]{\smash{{\SetFigFont{8}{9.6}{\updefault}{\color[rgb]{0,0,0}$v_{2a-2}$}%
}}}}
\put(3286,-16){\makebox(0,0)[rb]{\smash{{\SetFigFont{8}{9.6}{\updefault}{\color[rgb]{0,0,0}$v_{2a-1}$}%
}}}}
\put(1126,209){\makebox(0,0)[b]{\smash{{\SetFigFont{8}{9.6}{\updefault}{\color[rgb]{0,0,0}$\Delta_1$}%
}}}}
\put(1576,209){\makebox(0,0)[b]{\smash{{\SetFigFont{8}{9.6}{\updefault}{\color[rgb]{0,0,0}$\Delta_2$}%
}}}}
\put(1846,479){\makebox(0,0)[b]{\smash{{\SetFigFont{8}{9.6}{\updefault}{\color[rgb]{0,0,0}$\ldots$}%
}}}}
\put(1891,209){\makebox(0,0)[b]{\smash{{\SetFigFont{8}{9.6}{\updefault}{\color[rgb]{0,0,0}$\cdots$}%
}}}}
\put(2206,209){\makebox(0,0)[b]{\smash{{\SetFigFont{8}{9.6}{\updefault}{\color[rgb]{0,0,0}$\Delta_{b-1}$}%
}}}}
\put(3121,364){\makebox(0,0)[lb]{\smash{{\SetFigFont{8}{9.6}{\updefault}{\color[rgb]{0,0,0}$v_{2a}$}%
}}}}
\end{picture}%

%% file: EPS/le-gd.pstex_t
\begin{picture}(0,0)%
\includegraphics{le-gd.pstex}%
\end{picture}%
\setlength{\unitlength}{4144sp}%
\begingroup\makeatletter\ifx\SetFigFont\undefined
\def\x#1#2#3#4#5#6#7\relax{\def\x{#1#2#3#4#5#6}}%
\expandafter\x\fmtname xxxxxx\relax \def\y{splain}%
\ifx\x\y   
\gdef\SetFigFont#1#2#3{%
  \ifnum #1<17\tiny\else \ifnum #1<20\small\else
  \ifnum #1<24\normalsize\else \ifnum #1<29\large\else
  \ifnum #1<34\Large\else \ifnum #1<41\LARGE\else
     \huge\fi\fi\fi\fi\fi\fi
  \csname #3\endcsname}%
\else
\gdef\SetFigFont#1#2#3{\begingroup
  \count@#1\relax \ifnum 25<\count@\count@25\fi
  \def\x{\endgroup\@setsize\SetFigFont{#2pt}}%
  \expandafter\x
    \csname \romannumeral\the\count@ pt\expandafter\endcsname
    \csname @\romannumeral\the\count@ pt\endcsname
  \csname #3\endcsname}%
\fi
\fi\endgroup
\begin{picture}(3915,1912)(1396,-947)
\put(4141,-376){\makebox(0,0)[b]{\smash{\SetFigFont{8}{9.6}{rm}{\color[rgb]{0,0,0}$\vdots$}%
}}}
\put(5311,254){\makebox(0,0)[lb]{\smash{\SetFigFont{8}{9.6}{rm}{\color[rgb]{0,0,0}$v_0 = v_{\text{min}}$}%
}}}
\put(5041,-16){\makebox(0,0)[lb]{\smash{\SetFigFont{8}{9.6}{rm}{\color[rgb]{0,0,0}$v_{1}$}%
}}}
\put(4861,-196){\makebox(0,0)[lb]{\smash{\SetFigFont{8}{9.6}{rm}{\color[rgb]{0,0,0}$v_{2}$}%
}}}
\put(4501,-556){\makebox(0,0)[lb]{\smash{\SetFigFont{8}{9.6}{rm}{\color[rgb]{0,0,0}$v_{n-8}$}%
}}}
\put(4141,-916){\makebox(0,0)[b]{\smash{\SetFigFont{8}{9.6}{rm}{\color[rgb]{0,0,0}$v_{n-9}$}%
}}}
\put(3286,-16){\makebox(0,0)[rb]{\smash{\SetFigFont{8}{9.6}{rm}{\color[rgb]{0,0,0}$v_{2n-13}$}%
}}}
\put(3466,-196){\makebox(0,0)[rb]{\smash{\SetFigFont{8}{9.6}{rm}{\color[rgb]{0,0,0}$v_{2n-14}$}%
}}}
\put(3061,164){\makebox(0,0)[rb]{\smash{\SetFigFont{8}{9.6}{rm}{\color[rgb]{0,0,0}$v_{2n-12}$}%
}}}
\put(3781,-556){\makebox(0,0)[rb]{\smash{\SetFigFont{8}{9.6}{rm}{\color[rgb]{0,0,0}$v_{n-10}$}%
}}}
\put(2611,-286){\makebox(0,0)[b]{\smash{\SetFigFont{8}{9.6}{rm}{\color[rgb]{0,0,0}$f$}%
}}}
\put(1846,164){\makebox(0,0)[b]{\smash{\SetFigFont{8}{9.6}{rm}{\color[rgb]{0,0,0}$v_{\text{start}}$}%
}}}
\put(1396,254){\makebox(0,0)[rb]{\smash{\SetFigFont{8}{9.6}{rm}{\color[rgb]{0,0,0}$v_{\max} = v_{2n-5}$}%
}}}
\end{picture}